%% file: robustgraph_rev.tex
\documentclass[final]{elsart}

\usepackage{amsmath} 
\usepackage{amssymb}  
\usepackage{graphics}
\usepackage{color}
\usepackage{epsfig}
\usepackage{algpseudocode}
\usepackage{tikz,pgfplots}

\newcommand{\beq}{\begin{equation}}
\newcommand{\eeq}{\end{equation}}

\newcommand{\R}{\mathbb R}

\newcommand{\one}{\mathbf{1}}

\def\eps{\varepsilon}

\newcommand{\tF}{F}
\newcommand{\tG}{G}

\newcommand{\epstark}{\eps^{\star}}
\newcommand{\diagvec}{\mathrm{diagvec}}
\newcommand{\diag}{\mathrm{diag}}

\newcommand\Lap{L}
\newcommand\sym{\mathrm{Sym}}

\newcommand{\opt}{\mathrm{opt}}

\def\bcl{\color{black}}
\def\ecl{\color{black}}
\def\bde{\color{black}}
\def\ede{\color{black}}

\newtheorem{example}{Example}[section]

\newtheorem{theorem}{Theorem}[section]
\newtheorem{lemma}{Lemma}[section]

\newtheorem{remark}{Remark}[section]
\newtheorem{assumption}{Assumption}[section]

\begin{document}

\begin{frontmatter}

\title{Measuring the stability of spectral clustering}

\author[And]{Eleonora Andreotti}
\author[Ede]{Dominik Edelmann}
\author[Gug]{Nicola Guglielmi}
\author[Lub]{Christian Lubich}

\address[And]{Dipartimento di Fisica,
        Universit\`a degli Studi di Torino, Via Pietro Giuria 1,
        10125 Torino,  Italy
        {\tt\small andreotti@unito.it}}%
\address[Ede]{Mathematisches Institut,
       Universit\"at T\"ubingen,
       Auf der Morgenstelle 10,
       D--72076 T\"ubingen,
       Germany 
			 {\tt\small dominik.edelmann@na.uni-tuebingen.de}} 				
\address[Gug]{Gran Sasso Science Institute, via Crispi 7, 
I-67100 L'Aquila (Italy)
        {\tt\small nicola.guglielmi@gssi.it}}
\address[Lub]{Mathematisches Institut,
       Universit\"at T\"ubingen,
       Auf der Morgenstelle 10,
       D--72076 T\"ubingen,
       Germany 
			 {\tt\small lubich@na.uni-tuebingen.de}}

%

\begin{abstract}

As an indicator of the stability of spectral clustering of an undirected weighted graph into $k$ clusters, the $k$th spectral gap of the graph Laplacian is often considered. The $k$th spectral gap is characterized in this paper as an {\it unstructured} distance to ambiguity, namely as the minimal distance of the Laplacian to arbitrary symmetric matrices with vanishing $k$th spectral gap. As a conceptually more appropriate measure of stability, the {\it structured} distance to ambiguity of the $k$-clustering is introduced as the minimal distance of the  Laplacian to Laplacians of graphs with the same vertices and edges but with weights that are perturbed such that the $k$th spectral gap vanishes. To compute a solution to this matrix nearness problem, a two-level iterative algorithm is proposed that  uses a constrained gradient system of matrix differential equations in the inner iteration and a one-dimensional optimization of the perturbation size in the outer iteration.
	The structured and unstructured distances to ambiguity are compared on some example graphs. The numerical experiments show, in particular, that selecting the number $k$ of clusters according to the criterion of maximal stability can lead to different results for the structured and unstructured stability indicators.


\end{abstract}



\end{frontmatter}



\section{Introduction}\label{sec1}
Clustering --- identifying groups with similarities in a
large data set --- is a fundamental task in data analysis. 
In order to cluster an undirected weighted graph, several algorithms are considered in the literature.
Spectral clustering (originating with Fiedler \cite{fiedler}; see \cite{Sp07} 
and \cite{VL07} for more recent accounts)
offers significant  advantages over more traditional combinatorial techniques such as
direct $k$-means or single linkage. Quoting \cite{VL07}, ``{\it Results obtained by spectral clustering often outperform the 
	traditional approaches; spectral clustering is very simple to implement and can be solved efficiently 
	by linear algebra methods.}''

The stability of spectral clustering algorithms is heuristically --- and very conveniently --- associated with the spectral gap in the
graph Laplacian, that is the difference  between the 
$k$th and $(k+1)$-st smallest eigenvalues of the Laplacian matrix in the case of clusters built from $k$ eigenvectors. \bcl If these two eigenvalues coalesce, then the spectral $k$-clustering becomes ambiguous, and arbitrarily small perturbations may yield different clusterings.

In this paper we measure the stability of spectral $k$-clustering, for any number $k$ of clusters, by the 
{\it structured distance to ambiguity}, which is introduced here as the
distance (w.r.t. the Frobenius norm) of
the given Laplacian matrix to the set of 
Laplacian matrices that correspond to a weighted graph with the same vertices and edges but perturbed weights such that
the $k$th and $(k+1)$-st smallest eigenvalues of the graph 
Laplacian coalesce --- and for which therefore spectral clustering with $k$ clusters becomes ambiguous.
\ecl
On the other hand, the $k$th spectral gap divided by $\sqrt{2}$ \bcl will be \ecl characterized as the  distance of the Laplacian matrix to the much larger set of arbitrary, unstructured symmetric matrices with coalescing $k$th and $(k+1)$-st eigenvalues.

A stability indicator  is useful in selecting an appropriate number $k$ of clusters, by choosing $k$ such that the stability indicator is maximized. The structured distance to ambiguity, whose computation will be discussed in this paper, may be too costly computationally to be used as a general-purpose stability indicator in comparison with the spectral gap, which is available without further computations in spectral clustering. Nevertheless, comparing the structured distance to ambiguity and the spectral gap for representative examples within a class of graphs of interest gives insight into the reliability (or otherwise) of the spectral gap as a stability indicator.

The structured distance to ambiguity considered in this paper is similar in spirit to other structured robustness measures that arise in matrix and control theory, such as  stability radii, 
passivity distances,  distance to singularity, etc., \bde see \ede \cite{HinrichsenPritchard86,Higham88,Higham02,BennerVoigt14,GKL15,GLM17,KressnerVoigt15}.
The notion of stability considered here is, however, different from \cite{Luxburg10}, where a statistical perspective on clustering stability is developed.

A main objective of this paper is to show how the structured distance to ambiguity can be computed. 
The proposed algorithm is an iterative algorithm, where in each step a pair of adjacent eigenvalues and associated eigenvectors of the Laplacian of a graph with perturbed weights  are computed. For a large sparse graph (where the number of edges leaving any vertex is moderately bounded), these computations can typically be done with a complexity that is linear in the number of vertices. 

A feature of this algorithm in common with recent algorithms for eigenvalue optimization as given in \cite{GL11,GKL15,GLM17,GL17,AEGL19} is a two-level procedure for matrix nearness problems, where in an inner iteration a gradient flow drives perturbations to the original matrix of a {\it fixed} size into a (local) minimum of a nonnegative functional that depends on eigenvalues and eigenvectors, and in an outer iteration the perturbation size is optimized such that the functional becomes zero. 

The paper is organized as follows.
In Section~\ref{sec:problem} we recall basic concepts of spectral clustering. 
In Section~\ref{sec:robustness measure} we introduce the structured distance to ambiguity (SDA) as a measure of stability for spectral $k$-clustering and we characterize the $k$th spectral gap as an unstructured distance to ambiguity. 
In Section~\ref{sec:two-level method} we describe a numerical method to compute the SDA, which requires the solution of a structured matrix nearness problem. We propose a two-level iterative method, which is based on a gradient system of matrix differential equations in an inner iteration and a one-dimensional optimization of the perturbation size in an outer iteration, and we discuss algorithmic aspects. 
In Section~\ref{sec:numexp} we present the results of numerical experiments where the spectral gap (i.e., the unstructured distance to ambiguity) and the theoretically more appropriate measure of stability discussed here (i.e., the structured distance to ambiguity) are compared on some classes of graphs. We give a conclusion in Section~\ref{sec:conclusion}. 

\section{Spectral clustering}
\label{sec:problem}
Let $\mathcal{G}=(\mathcal{V},\mathcal E)$ be an undirected graph with vertex set $\mathcal{V}=\{1,\dots,n \}$ and 
edge set $\mathcal{E}\subset \mathcal{V}\times\mathcal{V}$.
We  assume that the graph is weighted, that is, a non-negative weight $w_{ij} =w_{ji} \ge 0$ is associated with each edge $(i,j)\in\mathcal{E}$ 
between two vertices $i$ and $j$.
We set $w_{ij}=0$ for $(i,j)\notin\mathcal{E}$.
The weighted adjacency matrix of the graph is the matrix 
$$
W=(w_{ij}) \in \R^{n\times n}.
$$
%
The degrees  $d_i = \sum_{j=1}^n w_{ij} $ are the elements of the diagonal matrix 
$$
D = \diag(d_i), \quad d_i=(W \one)_i, \qquad\hbox{where $\one:=(1,\ldots,1)^T \in \R^n$.}
$$
%
%
The (unnormalized) {\em Laplacian} of the graph is given by the matrix 
$L = \Lap(W)$, 
$$
L = D-W, \quad\text{ i.e.,}\quad L(W)= \diag (W  \one) - W .
$$
The Laplacian matrix $L$ is symmetric and positive semi-definite; since - by construction - 
$L\one=0$, $\lambda_1=0$ is the smallest  eigenvalue of $L$. 
Note that the matrix $L$ does not depend on (possible) diagonal elements of the matrix
$W$, which means that self-edges do not change the graph Laplacian.
The graph Laplacian and its eigenvectors provide important instruments to spectral
clustering, as stated by the following theorems.
\begin{theorem}[Bi-partition; Fiedler \cite{fiedler}] \label{thm:fiedler}
	Let $W \in \R^{n \times n}$ be the weight matrix of an undirected graph and $L$ the corresponding Laplacian matrix. 
	Let $0 = \lambda_1 \le \lambda_2 \le \ldots \le \lambda_n$ be the eigenvalues of $L$. Then, the graph is disconnected 
	if and only if $\lambda_2 = 0$. Moreover, if $0=\lambda_2<\lambda_3$, then the entries of the corresponding eigenvector 
	orthogonal to $\one$ assume only two different values, of different sign, which mark the membership to the two connected 
	components.
\end{theorem}

If $\lambda_2$ is a simple eigenvalue, then the corresponding eigenvector is known as the {\it Fiedler vector}. 
In spectral graph theory, inspired by Theorem \ref{thm:fiedler}, it is common to compute the second smallest
eigenvalue of $L$ and label the positive components of the Fiedler vector as belonging to one subset and the negative
ones to another subset, and in this way obtaining a natural partition of the graph.
However, this becomes unreliable when a small perturbation of the weights yields a coalescence of the eigenvalues 
$\lambda_2$ and $\lambda_3$. 

More generally we have the following result  (see e.g. \cite{VL07}). For a subset of vertices $C \subset \mathcal{V}$, we here denote the indicator vector $\one_C$ as the vector whose $\ell$th entry is equal   
to $1$ if $v_\ell \in C$ and is equal to zero otherwise. 

\begin{theorem}[$k$-partition] \label{thm:fiedlerk}
	Let $W \in \R^{n \times n}$ be the weight matrix of an undirected graph and $L$ the corresponding Laplacian matrix. 
	Then the multiplicity $k$ of the eigenvalue $0$ (the dimension of $\ker(L)$) equals the number of
	connected components $C_1, \ldots, C_k$ in the graph.
	The eigenspace of the eigenvalue $0$ is spanned by the indicator vectors $\one_{C_1},\ldots,\one_{C_k}$.
\end{theorem}

Nonempty sets $C_1,\ldots,C_k$ form a clustering of the graph if 
\[
C_i \cap C_j = \emptyset \quad \mbox{for} \ i,j=1,\ldots,k,\; i \neq j \quad \mbox{and} \quad \bigcup\limits_{i=1}^{k} C_i = \mathcal{V}.
\]
Similar to the case of two clusters,  this result motivates an algorithm for clustering a graph into $k$
components, which is reviewed in Algorithm~\ref{alg_sca}; see \cite{VL07}.
For the classical $k$-means algorithm we refer the reader e.g. to \cite{K14}.
Analogous results and algorithms are extended to the normalized Laplacian
\[
L_{\sym} := D^{-1/2}  L D^{-1/2} = I - D^{-1/2} W D^{-1/2}. 
\]


\renewcommand{\algorithmicrequire}{\textbf{Input:}}
\renewcommand{\algorithmicensure}{\textbf{Output:}}
\begin{algorithm}
	\caption{Unnormalized spectral clustering algorithm}\label{alg_sca}
	\algsetblock{Begin}{}{5}{}
	\begin{algorithmic}[1]
		\Require{weight matrix $W$, number $k$ of clusters}
		\Ensure{Clusters $C_1,\ldots,C_k$}
		\Begin
			\State Compute the Laplacian matrix $L=D-W$
			\State \bde Compute the eigenvectors $x_1,\ldots,x_k$ associated with the $k$ smallest nonzero eigenvalues  of $L$ \ede
			\State \bde Set $X = \left[ x_1 \ | \ x_2 \ | \ \ldots \ | \ x_k \right]$\ede
			\State \bde For $i=1,\ldots,n$,  define $r_i \in \R^{k}$ the vector given by the $i$-th row of~$X$ \ede
			\State Cluster the points $\left( r_i \right)_{i=1,\ldots,n}$ by the $k$-means algorithm into $k$ clusters $C_1,\ldots,C_k$
			\State Return $C_1,\ldots,C_k$ 
	\end{algorithmic}
\end{algorithm}

\section{Structured distance to ambiguity as a stability measure for spectral clustering}
\label{sec:robustness measure}
In order to evaluate the robustness of the spectral clustering algorithm it is essential to quantify the sensitivity to perturbations of the
invariant subspace associated with the $k$ eigenvectors used by the algorithm. Suppose that a small perturbation to the weights of the graph makes the eigenvalues $\lambda_k$ and $\lambda_{k+1}$ of the Laplacian coalesce. Then the eigenvector associated with the computed $k$th smallest eigenvalue can change completely and hence can yield a different clustering.

It has been suggested in the literature to use the spectral gap $\lambda_{k+1}-\lambda_k$ as an indicator of the stability of the $k$-clustering; cf.~\cite{VL07}.
For {\it unstructured} perturbations of the graph Laplacian, this is usually motivated by the Davis--Kahan theorem (see, e.g., \cite{SS90} and \cite{GVL13}), which tells us that the 
distance between the eigenspaces of the Laplacian of a graph and of any perturbed symmetric matrix has 
a bound that is proportional to the perturbation size and inversely proportional to the
spectral gap $\lambda_{k+1} - \lambda_{k}$ of the Laplacian $\Lap(W)$. 

In another direction, which is more related to the concepts studied here,
 the spectral gap can be characterized as an {\em unstructured} distance to ambiguity (up to the scaling factor $\sqrt2$) by the following result, where $\| A \|_F = \sqrt{\sum_{i,j} a_{ij}^2}$ denotes the Frobenius norm of a matrix $A=(a_{ij})$.

\begin{theorem} \label{thm:sg}
The $k$th spectral gap, divided by $\sqrt2$, is characterized as
\bde
\begin{equation}
\frac{\lambda_{k+1} - \lambda_{k}}{\sqrt{2}} = \min \left\{  \| \Lap(W)-\widehat L\|_F :\, \widehat{L} \in \mathrm{Sym}(n),\, \lambda_k(\widehat{L}) = \lambda_{k+1}(\widehat{L}) \right\}\,,
\label{dist-ambiguity-unstructured}
\end{equation}
where $\mathrm{Sym}(n)$ denotes the set of all symmetric real $n \times n$-matrices.\ede
\end{theorem}

{\it Proof.} Let here $L=\Lap(W)$ for brevity. \\
(a) We diagonalize $L=Q\Lambda Q^T$ with the diagonal matrix $\Lambda=\text{diag}(\lambda_i)$ of eigenvalues and the orthogonal matrix $Q$ formed of the eigenvectors. If we set
$\widetilde L = Q \widetilde \Lambda Q^T$ with $\widetilde\Lambda=\text{diag}(\widetilde\lambda_i)$ and with $\widetilde \lambda_i = \lambda_i$ for $i\notin \{k,k+1\}$  and coalescing $k$th and $(k+1)$-st eigenvalues $\widetilde \lambda_k=\widetilde \lambda_{k+1}=\tfrac12(\lambda_k+\lambda_{k+1})$,
then we obtain
\begin{align*}
\| L-\widetilde L \|_F^2 &= \| \Lambda - \widetilde \Lambda \|_F^2 
= (\lambda_k - \widetilde \lambda_k)^2 + (\lambda_{k+1} - \widetilde \lambda_{k+1})^2 
\\
&= 2\Bigl(\frac{\lambda_{k+1} - \lambda_{k}}2\Bigr)^2 = \tfrac12 (\lambda_{k+1} - \lambda_{k})^2
\end{align*}
so that
$$
\frac{\lambda_{k+1} - \lambda_{k}}{\sqrt{2}} =  \| L-\widetilde L\|_F .
$$
(b) Let now $\widehat L$ be an arbitrary symmetric matrix with ordered eigenvalues $\widehat\lambda_i$ such that the $k$th and $(k+1)$-st eigenvalues coalesce: $\widehat\lambda_k=\widehat\lambda_{k+1}$. We then have
$$
 \tfrac12 (\lambda_{k+1} - \lambda_{k})^2 = \min_\lambda \Bigl( (\lambda_k-\lambda)^2 + (\lambda_{k+1}-\lambda)^2\Bigr)
 \le   (\lambda_k-\widehat\lambda_k)^2 + (\lambda_{k+1}-\widehat\lambda_{k+1})^2 .
$$
By the Hoffman--Wielandt theorem \cite[p.\,368]{HoJ85},
$$
\sum_i (\lambda_i-\widehat \lambda_i)^2 \le \| L - \widehat L \|_F^2,
$$
and so we find
$$
\tfrac12 (\lambda_{k+1} - \lambda_{k})^2 \le \| L - \widehat L \|_F^2,
$$
which concludes the proof. \qed

The interpretation is that applying Algorithm \ref{alg_sca} is justified when the gap is relatively large, otherwise a small perturbation 
of the weight matrix may yield coalescence of $\lambda_k$ and $\lambda_{k+1}$ and significantly change the clustering.

The conceptual weakness of this  approach is that it considers {\it unstructured} perturbations of the Laplacian, as opposed to the {\it admissible} perturbations that are Laplacians $\Lap(\widehat W)$ of a perturbed weight matrix that preserves the symmetry and non-negativity, i.e., $\widehat w_{ij} = \widehat w_{ji}\ge 0$, and the sparsity pattern of~$W$, i.e., $\widehat w_{ij}=0$ for $(i,j)\notin \mathcal{E}$. In \eqref{dist-ambiguity-unstructured}, the minimizer is given by 
$$
\widetilde L = \Lap(W) +
\tfrac12 (\lambda_{k+1} - \lambda_{k}) x_k x_k^T - \tfrac12 (\lambda_{k+1} - \lambda_{k}) x_{k+1} x_{k+1}^T,
$$
where $x_k$ and $x_{k+1}$ denote normalized eigenvectors of $\Lap(W)$ to the eigenvalues $\lambda_k$ and $\lambda_{k+1}$, respectively.
Apart from exceptional cases, $\widetilde L$ is not the Laplacian of a graph.

We therefore propose the following stability measure for $k$-clustering, which is given by a {\em structured distance to ambiguity}, 
\bde
\begin{align*}
\begin{aligned}
	\delta_k(W) &= \min \big\{ && \| \Lap(W) - \Lap(\widehat{W}) \|_F :\, \widehat{W} \in \R^{n \times n} ,\, \lambda_{k}(\Lap(\widehat{W})) = \lambda_{k+1}(\Lap(\widehat{W})) \,,\\	
	& &&\widehat{w}_{ij} = \widehat{w}_{ji} \ge 0 \text{ for }(i,j) \in \mathcal{E} \text{ and } \widehat{w}_{ij} = 0 \text{ for }(i,j) \notin \mathcal{E} \big\} \,.
\end{aligned}
\end{align*}
\ede
%
\bcl Here, the Laplacian $L(W)$ is not compared with arbitrary symmetric matrices with vanishing $k$th spectral gap as in \eqref{dist-ambiguity-unstructured}, but with the Laplacians of weight matrices on the same graph with vanishing $k$th spectral gap. \ecl
In view of \eqref{dist-ambiguity-unstructured}, it is clear that 
$$
\delta_k(W) \ge (\lambda_{k+1}-\lambda_{k})/\sqrt2,
$$ and often $\delta_k(W)$ is substantially larger than the scaled spectral gap. In view of Theorem~\ref{thm:fiedler}, we further note that $\delta_1(W)$ is the Frobenius-norm distance of the Laplacian $\Lap(W)$ to that of a nearest disconnected graph.

For clustering a graph it is usually not known beforehand what the best choice of the number $k$ of clusters is. Asking for the most stable clustering of a graph (most stable with respect to admissible perturbations of the weights) determines the optimal number $k$ as
$$
k_{\rm opt}(W) = \arg\max_k \delta_k(W),
$$
where $k$ can be limited to  $k\le k_{\max}$ and $k_{\max}$ is given {\it a priori} or chosen such that $\lambda_{k_{\max}}$ is smaller than some threshold.
This criterion for the selection of $k$ is considered here instead of choosing $k$ such that the spectral gap $\lambda_{k+1} - \lambda_{k}$ is maximized. The latter is computationally cheaper, but appears conceptually less meaningful. In our numerical experiments in Section~\ref{sec:numexp} we will compare the two criteria for some families of graphs.

\begin{remark}
	There are obvious alternatives to the above stability measure:
	\begin{itemize}
		\item Instead of the unnormalized Laplacian $\Lap(W)=D-W$ with $D=\diag(W\one)$ we might work with the normalized Laplacian $I-D^{-1/2}WD^{-1/2}$.
		\item Instead of minimizing the perturbation of the Laplacian we might  minimize the perturbation of the weights, $\| W-\widehat W \|_F$.
		\item Instead of allowing for perturbations in every entry of the weight matrix $W$, perturbations might be restricted to selected entries.
	\end{itemize}
	In the following we concentrate on $\delta_k(W)$ as given above, but our algorithm for computing $\delta_k(W)$ is readily extended to the other cases.
\end{remark}

\section{Computation of the structured distance to ambiguity}
\label{sec:two-level method}
In this section we describe an approach to compute the structured distance to ambiguity $\delta_k(W)$ defined in the previous section.
\subsection{Outline of the computational approach}

Our approach is summarized by the following two-level method:
\begin{itemize}
	\item {\em Inner iteration:\/} Given $\eps>0$, we aim to compute a symmetric matrix $E=(e_{ij})\in\R^{n\times n}$ with the same sparsity pattern as $W$ (i.e., $e_{ij}= 0$ if $w_{ij}=0$), with $\Lap(E)=\diag(E\one)-E$ of unit Frobenius norm,  with $W+\eps E\ge 0$ (with entrywise inequality)  such that the difference between the {$(k+1)$-st} and $k$th smallest eigenvalues of $\Lap(W+\eps E)$ is minimal. The obtained minimizer is denoted by $E_{\eps}$.
	
	\item {\em Outer iteration:\/} We compute the smallest value of $\eps$ such that the $k$th and $(k+1)$-st eigenvalues of  $\Lap(W+\eps E_{\eps})$ coalesce.
\end{itemize}

In order to compute $E_{\eps}$ for a given $\eps>0$, we make use of a constrained gradient system for the functional
\begin{equation}\label{F-eps}
\tF_\eps(E) = \lambda_{k+1}  \bigl( \Lap(W+\eps E) \bigr) - \lambda_{k} \bigl( \Lap(W+\eps E) \bigr),
\end{equation}
under the constraints of unit Frobenius norm of $\Lap(E)=\diag(E\one)-E$, the non-negativity $W+\eps E\ge 0$ and symmetry and the sparsity pattern of~$E$. 

In the outer iteration we compute the smallest $\eps$ with $F_\eps(E_\eps)=0$, denoted~$\epstark$, by a combined Newton-bisection method.  The algorithm computes the perturbed weight matrix $W^\star=W+\epstark E_{\epstark}$ whose Laplacian has coalescent $k$th and $(k+1)$-st eigenvalues.
We then have
$$
\delta_k(W) = \epstark= \| \Lap(W)-\Lap(W^\star) \|_F.
$$
\bcl
\begin{remark} The constrained minimization of $F_\eps$ for a given $\eps$ is a nonconvex and nonsmooth optimization problem, which may have several local minima. By following a trajectory of the constrained gradient system into a stationary point, as we do here, we may thus end up in a local (instead of global) minimizer $E_\eps$. In our numerical experiments this situation was observed only in contrived examples \cite[Example~2]{AEGL19}. Using several trajectories starting from different random initial values can mitigate this potential problem. In any case, the algorithm computes an $\epstark= \| \Lap(W)-\Lap(W^\star) \|_F$ that is an upper bound to $\delta_k(W)$, and the bound can be made tighter by running more trajectories. 
\end{remark}
\ecl

\subsection{Gradient of the functional $\tF_\eps$}
\label{subsec:gradient}

We denote by $\langle X,Y \rangle = \mathrm{trace}(X^T Y)$ the inner product of matrices that corresponds to
 the Frobenius norm $\|X\|_F=\langle X,X \rangle^{1/2}$ on $\R^{n\times n}$.

For the edge set $\mathcal{E}$ and  for an arbitrary matrix $A =(a_{ij})\in \R^{n\times n}$, 
we define the symmetric projection onto the sparsity pattern given by $\mathcal{E}$ as
\[ 
P_\mathcal{E}(A)\big|_{ij} := 
\begin{cases}
\displaystyle{\tfrac12} \left( a_{ij} + a_{ji} \right) \, &\text{if } (i,j) \in \mathcal{E} \\[1mm]
0\, &\text{otherwise}.
\end{cases}
\]
We denote by $\sym(\mathcal{E})$ the space of symmetric matrices with sparsity pattern~$\mathcal{E}$ and note that $P_\mathcal{E}$ is the orthogonal projection from $\R^{n\times n}$ to $\sym(\mathcal{E})$: 
$$
\langle P_\mathcal{E}(A),W \rangle =\langle A,W \rangle \quad\text{ for all }\ 
W\in \sym(\mathcal{E}).
$$
The Laplacian operator is a linear map
$$
\Lap: \sym(\mathcal{E})\to \R^{n\times n}.
$$
Its adjoint with respect to the Frobenius inner product,
$$
\Lap^*: \R^{n\times n}\to \sym(\mathcal{E}),
$$
is given in the following lemma.
\begin{lemma}
	For $V\in \R^{n\times n}$,  let $\Lap^*(V)\in \sym(\mathcal{E})$  be defined by  $\langle \Lap^*(V),W \rangle = \langle V, L(W) \rangle$ for all  $W\in \sym(\mathcal{E})$. Then,
	$$
	\Lap^*(V) = P_\mathcal{E}(\diagvec(V)\one^T - V),
	$$
	where $\diagvec(V)\in \R^n$ is the vector of the diagonal entries of $V$. Furthermore, if $\mathcal{E}$ contains no auto-loops, i.e., $(i,i)\notin\mathcal{E}$ for all $i=1,\dots,n$, then we have for $W\in \sym(\mathcal{E})$,
	$$
	\Lap^*(\Lap(W)) = P_{\mathcal{E}} (d\one^T)  + W \quad\text{with }\  d=W\one.
	$$
\end{lemma}

{\bf Proof.}
	We have, for all $V\in \R^{n\times n}$ and $W\in \sym(\mathcal{E})$,
	\begin{align*}
	\langle V, L(W) \rangle &= \langle V, \diag(W\one)-W \rangle = \langle \diagvec(V), W\one \rangle - \langle V,W \rangle 
	\\
	&= \langle \diagvec(V)\one^T, W \rangle - \langle V,W \rangle = \langle P_\mathcal{E}(\diagvec(V)\one^T -V), W \rangle,
	\end{align*}
	\bcl which yields the formula for $\Lap^*(V)$. The formula for $\Lap^*(\Lap(W))$ then follows on inserting $V=L(W)=\diag(W\one)-W$ into this formula.\ecl
\qed

\bcl
The following lemma 
is basic for the derivation of the gradient system of matrix differential equations that we use to minimize the functional~$F_\eps$  of \eqref{F-eps}. Here we call a matrix $E\in \R^{n\times n}$ {\it $\eps$-feasible} if 
$$
E\in \sym(\mathcal{E}),\quad \| L(E) \|_F = 1, \quad W+\eps E \ge 0.
$$

\begin{lemma}\label{lem:lambda-dot} \cite[Lemma 3.2]{AEGL19} Let $\eps>0$ and
	consider a differentiable path of $\eps$-feasible matrices $E(t)$ $(t_0< t < t_1)$. Let $\lambda(t)$ be a continuous path of simple 
	eigenvalues of the Laplacian matrix $L(W+\eps E(t))$ with the associated eigenvector $x(t)$ of unit Euclidean norm.
	Then, $\lambda(\cdot)$ is differentiable and (denoting the derivative w.r.t. $t$ by a dot and omitting the argument $t$)
	\begin{equation}
	\label{G-eps}
	\dot\lambda = \langle  \Lap^*(xx^T) , \eps\dot{E} \rangle.
	\end{equation}
	We note that $\Lap^*(xx^T)=P_{\mathcal{E}}\left((x \bullet x) \one ^T  -x x^T \right)$,
	where $x \bullet x =(x_i^2)\in\R^n$ denotes the  vector of squares of the entries of the vector 
	$x=(x_i)\in\R^n$. 
\end{lemma}

An immediate consequence of Lemma \ref{lem:lambda-dot}  applied to the eigenvalues in the functional $\tF_\eps$ of \eqref{F-eps}  is the following result.

\begin{theorem}\label{thm:gradient} Let $\eps>0$.
Along a differentiable path $E(t)$ of $\eps$-feasible matrices such that the $k$th and $(k+1)$-st eigenvalues $\lambda_k(t)$ and $\lambda_{k+1}(t)$ of the Laplacian $L(W+\eps E(t))$ are simple eigenvalues, 
we have for the functional $\tF_\eps$ of \eqref{F-eps} that
$$
\frac1\eps \frac d{dt} F_\eps(E(t)) = \langle G_\eps(E(t)), \dot E(t) \rangle
$$
with the (rescaled) gradient
	\bde
	\begin{align}
	\label{eq:freegrad}
		\tG_\eps(E) &= \Lap^*(x_{k+1} x_{k+1}^T-x_k x_k^T) 
		\\ 
		&=P_{\mathcal{E}}\left( ( x_{k+1} \bullet x_{k+1} - x_k \bullet x_k ) \one ^T  - \left( x_{k+1} x_{k+1}^T - x_k x_k^T \right) \right),
		\nonumber
	\end{align}
\bcl	where $x_k$ and $x_{k+1}$ are the normalized eigenvectors associated with the eigenvalues $\lambda_{k}$ and $\lambda_{k+1}$, respectively, of $L(W+\eps E)$.\qed \ecl \ede
\end{theorem}

\bcl
The following result shows that the gradient $G_\eps(E)$ cannot be zero. This fact will be important later.

\begin{lemma}\label{lem:G-nonzero}
Let $x,y\in\R^n$ be normalized eigenvectors of a graph Laplacian $\widehat L =L(\widehat W)$ corresponding to different eigenvalues $\lambda \ne \mu$. Then,
$$
\Lap^*(xx^T -yy^T)\ne 0.
$$
\end{lemma}

{\bf Proof.} We show that assuming $\Lap^*(xx^T -yy^T) = 0$ leads to a contradiction. We would then have
$$
\tfrac12 (x_i^2+x_j^2) -x_ix_j - \bigl(\tfrac12(y_i^2+y_j^2)-y_iy_j\bigr)=0 \qquad \forall\, (i,j)\in\mathcal{E},
$$
i.e.,
\begin{equation} \label{G-zero}
\tfrac12(x_i-x_j)^2=\tfrac12(y_i-y_j)^2 \qquad \forall\, (i,j)\in\mathcal{E}.
\end{equation}
Multiplying with the weights $w_{ij}$ and summing over all $(i,j)\in\mathcal{E}$ yields, however,
$$
\tfrac12\sum_{(i,j)\in\mathcal{E}} w_{ij} (x_i-x_j)^2 = x^T \widehat Lx = x^T (\lambda x) = \lambda x^Tx = \lambda
$$
and in the same way
$$
\tfrac12\sum_{(i,j)\in\mathcal{E}} w_{ij} (y_i-y_j)^2 = \mu \ne \lambda,
$$
which is in contradiction to \eqref{G-zero}.
\qed
\ecl

\subsection{Gradient flow under the unit-norm constraint}

We consider now the constrained gradient system associated with the functional $\tF_\eps$ with  
the constraint that $\Lap(E(t))$ has unit Frobenius norm.  We then get the following 
matrix differential equation:
\begin{equation}\label{eq:ode}
\dot E = - \tG_\eps(E) + \kappa \,\Lap^*(\Lap(E)) \quad\hbox{ with }\quad
\kappa = {}\frac{\langle \tG_\eps(E),  \Lap^*(\Lap(E)) \rangle}{\|  \Lap^*(\Lap(E)) \|_F^2}.
\end{equation}
Here we note that the denominator cannot vanish, because $ \Lap^*(\Lap(E)) =0$ would imply that \bde $0=\langle \Lap^*(\Lap(E)), E \rangle = \langle \Lap(E), \Lap(E) \rangle = \| \Lap(E) \|_F^2 $ in contra\-diction to $\| \Lap(E) \|_F=1$. \ede 
The role of $\kappa$ is that of a Lagrange multiplier that ensures that the constraint $\| \Lap(E) \|_F=1$ is satisfied. The given formula for $\kappa$ is obtained as follows: Differentiating the constraint $\| \Lap(E) \|_F^2 = \langle \Lap(E),\Lap( E) \rangle =1$ gives 
$\langle \dot E,\Lap^*(\Lap(E)) \rangle=\langle \Lap(\dot E),\Lap( E) \rangle = \frac{1}{2} \frac{d}{dt} \| \Lap(E) \|_F^2 =0$. Taking the inner product of both sides of the differential equation with $\Lap^*(\Lap(E))$ yields
$$
0 = - \langle G_\eps(E), \Lap^*(\Lap(E)) \rangle + \kappa \, \| \Lap^*(\Lap(E)) \|_F^2,
$$
which gives $\kappa$.

\subsection{Non-negativity constraint}

It may happen that along the solution trajectory of \eqref{eq:ode} some positive entry of $W + \eps E$
becomes negative. In our experiments we never observed such a situation although this could potentially occur. In this subsection we explain how to deal with the further constraint $W + \eps E(t) \ge 0$ for all $t$.

One possibility would be to follow the lines of \cite{AEGL19} and consider KKT conditions by managing
active constraints and integrating a piecewise smooth system of differential equations.
Another possibility is to add a penalization term such as the following:
\[
Q_\eps(E) = \frac12 \sum\limits_{(i,j) \in \mathcal{E}} \left( w_{ij} + \eps e_{ij} \right)_{-} ^2
\]
where $(a)_{-} = \min\left( a , 0 \right)$, and to minimize the functional
\[
F_{\eps,c}(E) = \tF_\eps(E) + c\, Q_\eps(E)
\] 
for increasing $c$. We have
\bde
\[
\nabla Q_\eps(E) = {}-\eps \left( W + \eps E \right)_{-},
\]
\ede
giving an extra term to the gradient system \eqref{eq:ode}.

With the notation
\begin{equation} \label{eq:grad}
G_{\eps,c}(E) = \tG_\eps(E) \bcl -  c\left( W + \eps E \right)_{-},
\end{equation}
the differential equation \eqref{eq:ode} is then replaced by
\begin{equation}\label{eq:ode2}
\dot E = - G_{\eps,c}(E) + \kappa \,\Lap^*(\Lap(E)) \quad\hbox{ with }\quad
\kappa = {}\frac{\langle G_{\eps,c}(E),  \Lap^*(\Lap(E)) \rangle}{\|  \Lap^*(\Lap(E)) \|_F^2}.
\end{equation}

\subsection{Monotonicity of the spectral gap and stationary points}

The following monotonicity result follows directly from the construction of the gradient system.

\begin{theorem}
	Let $E(t)$ \bde with $L(E(t))$  of unit Frobenius norm \ede satisfy the differential equation \eqref{eq:ode2} with $G_{\eps,c}(E)$ of \eqref{eq:grad}. 
	Then, $F_{\eps,c}(E(t))$ 
	decreases monotonically with $t$.
	Moreover, if $W + \eps E(t) \ge 0$, then
	$$
	\frac{d}{dt} (\lambda_{k+1}(t) - \lambda_{k}(t) ) \le 0,
	$$
	where $\lambda_{k}(t)$ and $\lambda_{k+1}(t)$ are the $k$th and $(k+1)$-st eigenvalues of $\Lap(W+\eps E(t))$.
\end{theorem}

{\bf Proof.}
	Using \eqref{eq:ode2} we obtain, with $G=G_{\eps,c}(E)$,
	\bde
	\begin{equation}\label{F-monotone}
	\begin{aligned}
		\frac1\eps\, \frac{d}{dt} F_{\eps,c} \left( E(t) \right)  &= \langle G, \dot E \rangle 
		= \langle G, - G + \kappa \,\Lap^*(\Lap(E)) \rangle
		\\
		&= - \| G \|_F^2 + 
		\frac{\langle G, \Lap^*(\Lap(E)) \rangle^2}{\| \Lap^*(\Lap(E)) \|_F^2} \le 0,
	\end{aligned}
	\end{equation} \ede
	where the final inequality follows directly from the Cauchy--Schwarz inequality. This yields the monotonic decay of $F_{\eps,c}(E(t))$ and hence also the second statement, since $F_{\eps,c}(E(t)) = \lambda_{k+1}(t) - \lambda_{k}(t)$ if $W + \eps E(t) \ge 0$.
\qed

Stationary points of \eqref{eq:ode2} are characterized as follows.

\begin{theorem} \label{thm:stat}
	The following statements are equivalent along solutions of \eqref{eq:ode2}.
	\begin{enumerate}
		\item $\frac{d}{dt} F_{\eps,c} \left( E(t) \right) =0$.
		\item $\dot E =0$.
		\item \bcl $G_{\eps,c}(E)$ is a real multiple of $\Lap^*(\Lap(E))$. \ecl
	\end{enumerate}
\end{theorem}

{\bf Proof.} 
	By  the strict Cauchy--Schwarz inequality in \eqref{F-monotone}, $\frac{d}{dt} F_{\eps,c} \left( E(t) \right)$ can be zero only if $G$ is a multiple of $\Lap^*(\Lap(E))$, and then $\langle G,\dot E \rangle=0$ with $\dot E = -G +\kappa \Lap^*(\Lap(E))$ implies  $\dot E=0$. 
\qed

\bde
The inner iteration is formulated as pseudocode in Algorithm~\ref{alg_inner}, where a projected explicit Euler method is used for the discretization of the differential equation \eqref{eq:ode2}.
\begin{algorithm}[h!]
	\caption{Inner iteration}
	\label{alg_inner}
	\algsetblock{Begin}{}{3}{}
	\begin{algorithmic}[1]
		\Require{$W$, $k$, $\eps$, $c$ (penalization parameter), $E_0$ (initial perturbation), tol (tolerance), $h$ (step size for gradient flow)}
		\Ensure{$E_{\eps,c}$, $F_{\eps,c}(E_{\eps,c})$}
		\Begin
		\State Compute $\Lap(W+\eps E_0)$, $F_{\eps,c}(E_0)$
		\State Initialize $F^{(new)}_{\eps,c} = F_{\eps,c}(E_0)$, $F^{(old)}_{\eps,c} = F^{(new)}_{\eps,c}+1$, $E=E_0$
		\While{$|F^{(new)}-F^{(old)}|/|F^{(new)}|>\mathrm{tol}$}
			\State Compute $G_{\eps,c}(E)$, $\Lap^*(L(E))$ and $\kappa$
			\State Compute $\dot{E} = -G_{\eps,c}(E)+\kappa \Lap^*(L(E))$
			\State Update $E := E+h \dot{E}$
			\State Normalize $E$ such that $\| \Lap(E)\|_F = 1$
			\State Compute $\Lap(W+\eps E)$ and $F_{\eps,c}(E)$
			\State Update $F^{(old)}_{\eps,c} = F^{(new)}_{\eps,c}$, $F^{(new)}_{\eps,c} = F_{\eps,c}(E)$
		\EndWhile	
		\State Return $E_{\eps,c}=E$ and $F_{\eps,c}(E_{\eps,c})=F_{\eps,c}^{(new)}$
	\end{algorithmic}	
\end{algorithm}
\begin{remark} This algorithm can be refined by using an adaptive stepsize strategy as in \cite[Section~5.3]{AEGL19}.
	In particular, when the gradient flow is close to the minimum, there can be oscillations in the functional $F_{\eps,c}(E)$, which are due to the discrete steps in the gradient system. This can be avoided by repeating lines 7 to 9 in Algorithm~\ref{alg_inner} with a smaller stepsize $h$ if the functional $F_{\eps,c}(E)$ grows.
\end{remark}
\ede

\bcl
\subsection{One-dimensional optimization over $\eps$}
\label{sec:outer}

Let $E_{\eps,c}$ denote the minimizer of the functional $F_{\eps,c}$, and let $c$ be so large that $W+\eps E_{\eps,c}\ge 0$, which implies that $F_{\eps,c}(E_{\eps,c})=F_{\eps}(E_{\eps,c})=:f(\eps)$.
Let $\epstark>0$ be the minimal zero of $f(\eps)$.
Generically we expect that for a given perturbation size $\eps<\epstark$,  
the eigenvalues ${\lambda_{k} (W+\eps E_{\eps,c})> 0}$ and 
$\lambda_{k+1} (W+\eps E_{\eps,c}) >0$ are simple. 
%
If so, then $f_c(\eps):=F_{\eps,c}(E_{\eps,c})$ is a smooth function of $\eps$ and we 
can exploit its regularity to obtain a fast iterative method to converge to
$\epstark$ from the left.
Otherwise we can use a bisection technique  to approach $\epstark$.


\begin{assumption}
	We assume that the $k$th and $(k+1)$-st smallest eigenvalues  of the Laplacian 
	$\Lap(W + \eps E_{\eps,c})$  are \emph{simple}. Moreover,  $E_{\eps,c}$ is assumed 
	to be  a smooth function of $\eps$ in some interval. 
	\label{assumpt}
\end{assumption}

We then have the following result.

\begin{lemma}
	\label{lem:der}
	Under Assumption~{\rm \ref{assumpt}}, the function $f_c(\eps)=F_{\eps,c}(E_{\eps,c})$ 
	is differentiable and its derivative equals (with ${\phantom{a}'}= d/d\eps$)
	\begin{equation}
	 f_c'(\eps) = \langle G_{\eps,c}(E_{\eps,c}), E_{\eps,c} \rangle. 
	\label{eq:derFdeps} 
	\end{equation}
\end{lemma}
{\bf Proof.} 
	Differentiating $f_c(\eps)=F_{\eps,c}( E_{\eps,c})$ with respect to $\eps$ we obtain as in Theorem~\ref{thm:gradient}   and noting that $\frac{d}{d\eps} (\eps E_{\eps,c}) =
	E_{\eps,c}+ \eps E'_{\eps,c}$ with $E'_{\eps,c} = \frac{d}{d\eps} E_{\eps,c}$, and setting $G(\eps)=G_{\eps,c}(E_{\eps,c})$, that
	$$
	f_c'(\eps) =\langle G(\eps), E_{\eps,c}+ \eps E'_{\eps,c} \rangle.
	$$	
	Now we use the property of minimizers, stated by Theorem~\ref{thm:stat}, and conclude
	$$ 
	G(\eps)=\kappa(\eps) \Lap^*(\Lap(E_{\eps,c})),  \quad\text{ with } \  
	\kappa(\eps) = \frac{\langle G(\eps),  \Lap^*(\Lap(E_{\eps,c})) \rangle}{\|  \Lap^*(\Lap(E_{\eps,c})) \|_F^2},
	$$
	which yields that, with 
	\begin{eqnarray*}
	\langle G(\eps), E'_{\eps,c} \rangle & = & \kappa(\eps) \langle \Lap^*(\Lap(E_{\eps,c})), E'_{\eps,c} \rangle = \kappa (\eps) \langle \Lap(E_{\eps,c})), \Lap(E'_{\eps,c}) \rangle =
	\\
	& = &  \frac{\kappa(\eps)}2 \frac{d}{d\eps} \| \Lap(E_{\eps,c}) \|_F^2 =0,
	\end{eqnarray*}
	since $\Lap(E_{\eps,c})$ is of unit norm for all $\eps$. \bcl So we have $f_c'(\eps) =\langle G(\eps), E_{\eps,c} \rangle$. \ecl
\qed

\bcl
Using Theorem~\ref{thm:stat} and Lemmas~\ref{lem:G-nonzero} and~\ref{lem:der}, we obtain the following result.

\begin{theorem}\label{thm:eps-monotone} If Assumption~\ref{assumpt} is satisfied for $\eps\in(0,\epstark)$ and $c$ is sufficiently large so that $W+\eps E_{\eps,c}\ge 0$ for all $\eps$, then $E_\eps:=E_{\eps,c}$ is independent of~$c$ and $f(\eps): =
F_\eps(E_{\eps})= f_c(\eps)$ has the following properties:
\begin{enumerate}
\item $f(\eps)$ is positive and strictly monotonically decreasing for $\eps\in(0,\epstark)$, and
$$
f'(\eps) = - \frac{\| G_\eps(E_{\eps})\|_F}{\| L^*(L(E_{\eps})) \|_F}\,.
$$
\item $f(\eps)=0$ for $\eps\ge\epstark$.
\end{enumerate}
\end{theorem}

{\bf Proof.} 
The independence of $c$ is obvious if $W+\eps E_{\eps,c}\ge 0$, since then the penalty term $Q_\eps(E_{\eps,c})=0$.

(1) We consider first small $\eps$ near 0. For any $E$ with $\| L(E) \|_F =1$ we have
$$
F_\eps(E) = F_0 +\eps \langle G_0,E\rangle + O(\eps^2),
$$
where $F_0$ and $G_0$ are defined as $F_\eps$ and $G_\eps$ for $\eps=0$ (they are then independent of $E$).
With the choice
$$
E = E_0 := -\frac{G_0}{\| L(G_0)  \|_F}
$$
we clearly have $\langle G_0,E\rangle = - \| G_0 \|_F^2 / \| L(G_0)  \|_F < 0$. Since $f(\eps)=F_\eps(E_\eps) \le F_\eps(E_0)$ and $G(\eps)=G_\eps(E_\eps)=G_0 +O(\eps)$, we conclude that for sufficiently small $\eps>0$ we have (with $E(\eps)=E_\eps$)
\begin{equation}\label{GE-neg}
\langle G(\eps), E(\eps) \rangle < 0.
\end{equation}
We now show that this inequality actually holds also for all $\eps\in (0,\epstark)$. By Theorem~\ref{thm:stat} we have
$$
G(\eps) = \kappa(\eps) \Lap^*(\Lap(E(\eps))),
$$
where $\kappa(\eps)\ne 0$ by Lemma~\ref{lem:G-nonzero}. We note that this implies
$$
|\kappa(\eps)| = \frac{\| G(\eps)\|_F}{\| \Lap^*(\Lap(E(\eps))) \|_F}\,.
$$
We then have
$$
\langle G(\eps), E(\eps) \rangle = \kappa(\eps) \langle \Lap^*(\Lap(E(\eps))), E(\eps) \rangle = 
\kappa(\eps) \| \Lap(E(\eps)) \|_F^2 = \kappa(\eps) \ne 0.
$$
Hence $\langle G(\eps), E(\eps) \rangle$ does not change sign, and so we have \eqref{GE-neg} for all $\eps\in (0,\epstark)$. By Lemma~\ref{lem:der}, this further yields
$$
f'(\eps) = \langle G(\eps), E(\eps) \rangle = \kappa(\eps)  <0.
$$
This proves statement (1) of the theorem.

(2) 
Once we have reached the minimum value $\epstark$ such that $f(\epstark) = 0$ the function remains identically zero. In order to see this
recall that $\Lap(W+\epstark E^\star)$ (with $E^\star = E(\epstark)$ the extremizer) has coalescent eigenvalues $\lambda_k=\lambda_{k+1}$.
Then consider the perturbation, for arbitrary  $\delta> -1$,
\begin{equation} \label{delta}
\eps E = \epstark E^\star + \delta\, (W + \epstark E^\star) ,
\end{equation}
where $\eps>0$ is determined by the normalization condition $\| L(E) \|_F=1$.
Clearly, $\Lap(W + \eps E)=(1+\delta)\Lap(W+\epstark E^\star)$ still has coalescent eigenvalues $\lambda_k=\lambda_{k+1}$,
so that  we have $f(\eps)=F_\eps(E)=0$. This implies $\eps\ge \epstark$, since  $\epstark$ is the smallest 
solution $\eps$ of the problem.
Evidently when $\delta \rightarrow \infty$ then $\eps \rightarrow \infty$, so that  $f(\eps)$
vanishes on the whole half-line $\eps \ge \epstark$. 
\qed

\begin{remark} Forming the Laplacian and taking the square of the Frobenius norm on both sides of \eqref{delta}, noting the normalizations
$\| L(E) \|_F = \| L(E^\star) \|_F =1$ and introducing $W^\star = W+\epstark E^\star$, we have
$$
(\epstark)^2 \le \eps^2 = (\epstark)^2 + 2\delta\, \langle L(W^\star)-L(W), L(W^\star) \rangle + \delta^2 \| L(W^\star) \|_F^2.
$$
This inequality holds true for all $\delta>-1$ if and only if the orthogonality relation
\begin{equation}\label{orth-rel}
 \langle L(W^\star)-L(W), L(W^\star) \rangle =0
\end{equation}
is satisfied. This necessary condition for optimality provides a simple check for a numerically obtained extremizer.
\end{remark}
\ecl


\subsection{A Newton bisection method}
For $\eps < \epstark$, we make use of the standard Newton iteration
\begin{equation}
\hat\eps = \eps - \frac{f_c(\eps)}{f_c'(\eps)},
\label{eq:Newton}
\end{equation}
to get an updated value $\hat\eps$.
In a practical algorithm it is useful to couple the Newton iteration \eqref{eq:Newton} with
a bisection technique, in the same way as for the method presented in \cite{AEGL19}.
We adopt a tolerance tol which allows us to distinguish
whether $\eps < \eps^\star$, in which case we may use the derivative formula and
perform the Newton step, or $\eps > \eps^\star$, so that we have to make use of bisection.
The method is formulated in Algorithm \ref{alg_dist}.

\begin{algorithm} [ht!]
	\caption{Newton-bisection method}
	\label{alg_dist}
	\algsetblock{Begin}{}{3}{}
	\begin{algorithmic}[1]
	\Require{$W$, $m_{\max}$ (max number of iterations), 
		tol (tolerance),\; \bde $\eps_{\rm lb}$ and $\eps_{\rm ub}$ (starting values for the lower and upper
		bounds for $\eps^\star$), $\eps_0 \in [\eps_{lb},\eps_{ub}]$, $c$ (regularization parameter)}\ede
	\Ensure{$\eps^\star$ (upper bound for the distance), $E(\eps^\star)$}
	\Begin
		\State Compute $E(\eps_0)$ \bde by applying Algorithm~\ref{alg_inner} \ede\;
		\State Set $m=0$\;
		\While{$m \le m_{\max}$}
			\If{$f_c(\eps_m) < {\rm tol}$}
				\State Set $ \eps_{\rm ub} = \min(\eps_{\rm ub},\eps_m)$\;
				\State Set $\eps_{m+1} = (\eps_{\rm lb} + \eps_{\rm ub})/2$ \ (bisection step)
				\Else
				\State Set $ \eps_{\rm lb} = \max(\eps_{\rm lb},\eps_m)$\;
				\State Compute $f_c(\eps_m)$ and $f_c'({\eps_m})$\;
				\State Compute $\eps_{m+1} = \displaystyle{\eps_m - \frac{f_c(\eps_m)}{f_c'(\eps_m)}}$ \ (Newton step)
			\EndIf
			\If{$\eps_{m+1} \not\in (\eps_{\rm lb},\eps_{\rm ub})$}
				\State Set $\eps_{m+1} = (\eps_{\rm lb} + \eps_{\rm ub})/2$
			\EndIf
			\If{$m=m_{\max}$ {\rm \bf or} $\eps_{\rm ub}-\eps_{\rm lb} < {\rm tol}$}
				\State Return $\eps_{m+1}$ {\rm \bf and} the interval $[\eps_{\rm lb},\eps_{\rm ub}]$\;
				\State {\rm \bf Stop}
			\Else
				\State Set $m=m+1$
			\EndIf
			\State Compute $E(\eps_m)$ \bde by applying Algorithm~\ref{alg_inner} \ede\;
		\EndWhile
	\State Return $\eps^\star = \eps_m$\;
	\end{algorithmic}	
\end{algorithm}

The overall method is formulated in Algorithm \ref{alg_dist_c}. \bcl In our experience, usually  it already terminates
\bcl with $c_0=0$.\ecl

\begin{algorithm}[h!]
\bcl
	\caption{The overall algorithm}
	\label{alg_dist_c}
	\algsetblock{Begin}{}{6}{}
	\begin{algorithmic}[1]
	\Require{Matrix $W$ and
		 increasing values $c_0=0,c_1,\ldots,c_m$ of penalization parameters}
	\Ensure{$\epstark$ (upper bound for the distance), $E^\star=E(\epstark)$}
	\Begin
		\State Set $\ell=0$\;
		\While{$\ell \le m$}	
		\State Find $\eps_\ell = \min_{\eps} \{ \eps \,|\, F_{\eps,c_\ell} (E) =0  \}$ 
		by applying Algorithm~\ref{alg_dist}\; 
		\State Set $E_\ell$ the minimizer (normalized by $\| L(E_\ell) \|_F=1$)\; 
		\If{$W+\eps_\ell E_\ell \ge 0$}
			\State Return $\epstark=\eps_\ell$, $E^\star=E_\ell$\;
	                 \State {\rm \bf Stop}
		\EndIf
		\State Set $\ell = \ell +1$
	\EndWhile
		\end{algorithmic}
\ecl
\end{algorithm}

\pagebreak[3]

\subsection{Effective monotonicity with respect to $\eps$} 
\bcl
Assume we have integrated equation \eqref{eq:ode2} with $\eps=\eps_0$ and we increase to a new value
$\eps_1 > \eps_0$. It may not be possible to take  the minimizer $E_{\eps_0}$ (with
$\| L(E_{\eps_0}) \|_F=1$ and $W+\eps_0 E_{\eps_0}\ge 0$) computed in the previous step as the initial value for the constrained gradient flow \eqref{eq:ode2}, since we might no longer have  $W+\eps_1 E_{\eps_0}\ge 0$ for $\eps_1>\eps_0$.
We therefore first
integrate the free gradient system for $F_{\eps_1,c}$ with initial value proportional to  $E_{\eps_0}$,
$$
\dot E =  -\tG_{\eps_1,c}(E), \qquad
		E(0) = \displaystyle{\frac{\eps_0}{\eps_1}} E_{\eps_0},
$$
that is, the differential equation \eqref{eq:ode2} without the Lagrange multiplier $\kappa$ that enforces the
norm constraint. Note that $\| L(E(0)) \|_F =\eps_0/\eps_1< 1$. Along the solution $E(t)$, the functional $F_{\eps_1,c}(E(t))$ decreases by Theorem~\ref{thm:gradient}. For the norm of $L(E(t))$ we have
\begin{align*}
&\frac12\,\frac{d}{dt}\Big|_{t=0} \| L(E(t)) \|_F^2 = \langle L(\dot E(0)), L(E(0)) \rangle \\
&= \langle -\tG_{\eps_1,c}(E(0)), L^*(L(E(0))) \rangle
=  \langle -\tG_{\eps_0,c}(E_{\eps_0}), \frac{\eps_0}{\eps_1} L^*(L(E_{\eps_0})) \rangle. 
\end{align*}
Since $W+\eps_0 E_{\eps_0}\ge 0$, we have $\tG_{\eps_0,c}(E_{\eps_0})= \tG_{\eps_0}(E_{\eps_0})$, and from Theorem~\ref{thm:stat} and the proof of 
Theorem~\ref{thm:eps-monotone} we know that
$$
\tG_{\eps_0}(E_{\eps_0}) = \kappa(\eps_0) L^*(L(E_{\eps_0}))\quad\text{ with } \quad  \kappa(\eps_0)<0.
$$
Hence we obtain
$$
\frac{d}{dt}\Big|_{t=0} \| L(E(t)) \|_F^2 >0,
$$
and so the norm of $L(E(t))$ increases with $t$, at least initially. If $\eps_1$ is sufficiently close to $\eps_0$, we are thus guaranteed to reach
an instant $\bar t$ where 
\[
\| \Lap(E( \bar t )) \| =  1.
\] 
We then continue with the norm-preserving differential equation \eqref{eq:ode2},
$$
\dot E = - G_{\eps,c}(E) + \kappa \,\Lap^*(\Lap(E)) \quad\hbox{ with }\quad
\kappa = {}\frac{\langle G_{\eps,c}(E),  \Lap^*(\Lap(E)) \rangle}{\|  \Lap^*(\Lap(E)) \|_F^2}
$$
with the initial value $E( \bar t )$, which then has $\Lap(E( t ))$ of unit norm for all $t$.

In this way, alternating between  the free and the norm-constrained gradient flow, the computed functional $F_\eps( E(t) )$ is made decreasing with respect to both $t$ and $\eps$. \ecl

\subsection{Choice of the inital perturbation size $\eps_0$ and the initial perturbation matrix} 
While one might just start with a random perturbation, an educated guess starts from the normalized 
free gradient  $E^0=- G_\eps(0)/\| \Lap(G_\eps(0)) \|_F$.
We propose to choose as a first guess the scaled spectral gap $\eps_0=\bigl(\lambda_{k+1} \left(\Lap(W) \right) - \lambda_{k} \left(\Lap(W) \right)\bigr)/\sqrt2$, which underestimates 
the optimal $\epstark=\delta_k(W)$ and therefore leads us into the regime where Newton rather than bisection is used for 
updating $\eps$ in the outer iteration. 

%
%

\section{Numerical experiments}
\label{sec:numexp}
In this section we compare the behavior of the spectral gap and the structured distance to ambiguity (SDA) as stability indicators. First, we  determine the optimal numbers of clusters by the criterion of maximal stability with both stability indicators in a family of stochastic block models with varying edge probabilities, alongside a reduced model with similar behavior. Then we examine their performance on a model with randomly generated numbers near given centers, where the optimal number of clusters is given a priori but is here determined by the criterion of maximal stability with both stability indicators.

\subsection{Stochastic block model and reduced model}
The stochastic block model (SBM) is a model of generating random graphs that tend to have communities. It is an important model in a wide range of fields ranging from sociology to physics (see \cite{HLL83} and, e.g., \cite{abbe2015community,amini2018semidefinite}).

Here we consider a stochastic block model with the following parameters:
\begin{itemize}
	\item the number $n$ of vertices,
	\item a partition of the vertex set $\{1,\ldots,n\}$ into $r$ disjoint subsets $C_1,\ldots,C_r$ called communities,
	\item a symmetric  $(r\times r)$-matrix $P=(p_{ij})$  containing edge probabilities.
\end{itemize}
The graph is then sampled in the way that any two vertices $v_i \in C_i$ and $v_j \in C_j$ are connected with probability $p_{ij}$.

For the stability question considered in this paper, this model can be reduced to a model with $2 r$ vertices. Assume that the probability matrix $P$ of the SBM is such that $p_{ii}=1$ for all $1 \le i \le r$. In this case, any two vertices $v$ and $w$ belonging to the same community are connected. If $P = I$, the SBM thus generates a disconnected graph with $r$ communities, which are complete graphs of size $|C_1|,\ldots,|C_r|$. We represent this graph by a graph with vertex set $\{1,\ldots,2r\}$ such that the vertices $2k-1$ and $2k$ are connected with weight $|C_k|$. The edge probabilites $p_{ij}$ in case $P \ne I$ are then represented by inserting matrices $\mu_{ij} I_2$ in the respective part of the weight matrix in the reduced model, where $\mu_{ij}$ is an appropriate function of $p_{ij}$ that takes the community sizes into account.

\begin{figure}
	\centering
	\includegraphics[scale=0.55]{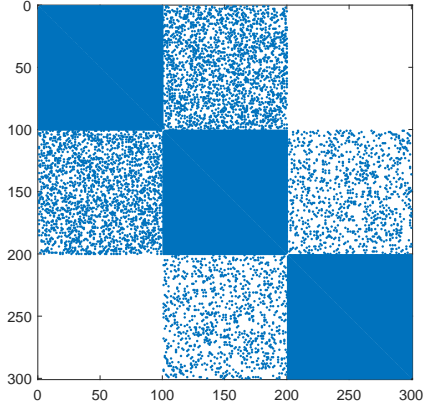}
	\caption{Adjacency matrix of a stochastic block model with $300$ vertices and three communities, sampled with the probability matrix \eqref{eq: probability matrix}.}\label{fig:spySBMmatrix}
\end{figure}

\begin{figure}
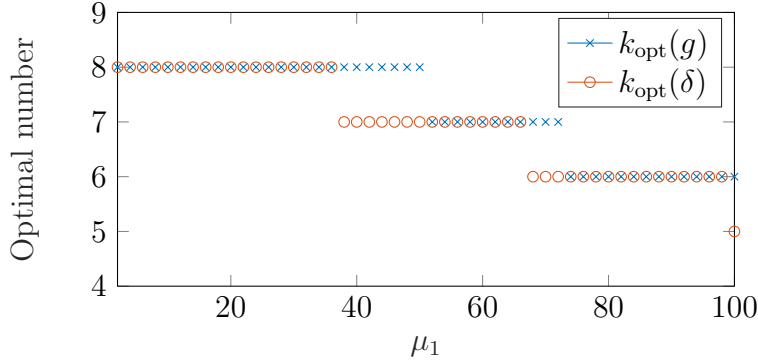

	\centering
	\include{reducedSBM1}
	\caption{Optimal number of clusters in Example~\ref{ex: chain sbm} provided by the SDA and by the spectral gaps.}\label{fig:kOpt}
\end{figure}

\begin{figure}
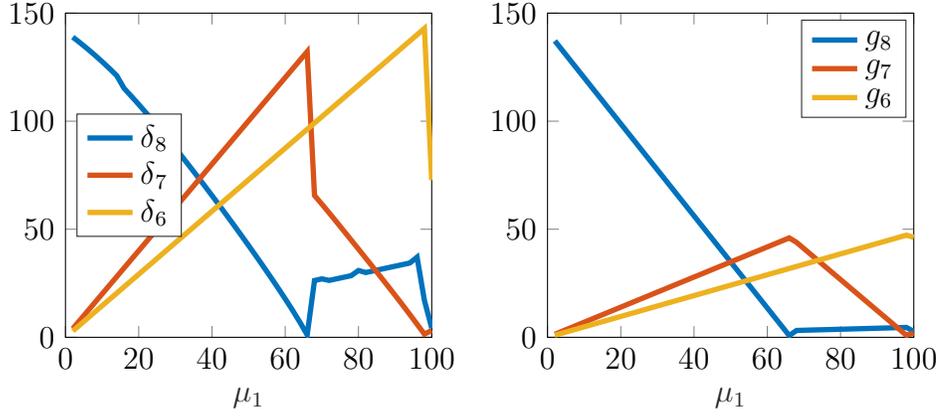

	\centering
	\include{reducedSBM2}
	\caption{Structured distances to ambiguity $\delta_6, \delta_7$ and $\delta_8$ (left) and spectral gaps $g_6, g_7$ and $g_8$ (right) in Example~\ref{ex: chain sbm} for the reduced model.}\label{fig:kVals}
\end{figure}

\begin{figure}
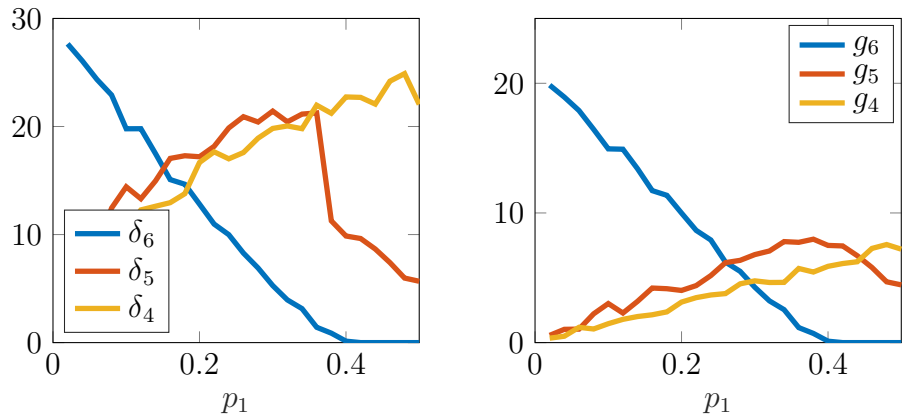

	\centering
	\include{SBM2}
	\caption{Structured distances to ambiguity $\delta_4, \delta_5$ and $\delta_6$ (left) and the spectral gaps $g_4, g_5$ and $g_6$ (right) in Example~\ref{ex: chain sbm} for the full SBM model.}\label{fig:kVals_SBM}
\end{figure}

To illustrate the construction of this reduced model, consider an SBM with $n=300$ vertices, three communities of size $100$ and the probability matrix
\begin{align}\label{eq: probability matrix}
	P = \begin{pmatrix} 1 & 0.2 & 0 \\
	0.2 & 1 & 0.1\\
	0 & 0.1 & 1\end{pmatrix}\,.
\end{align}
Figure \ref{fig:spySBMmatrix} shows the adjacency matrix of one sample of the corresponding SBM. This $(300 \times 300)$-matrix is represented by the $(6 \times 6)$-matrix
\begin{align}\label{eq: reduced block model matrix}
	W = \begin{pmatrix}
	0 & 100 & 20 & 0 &  &  \\ 
	100 & 0 & 0 & 20 &  &  \\ 
	20 & 0 & 0 & 100 & 10 & 0 \\ 
	0 & 20 & 100 & 0 & 0 & 10 \\ 
	&  & 10 & 0 & 0 & 100 \\ 
	&  & 0 & 10 & 100 & 0
	\end{pmatrix} \,.
\end{align}
With regards to clustering stability, we observe a similar behavior for the full and the reduced model, as the following example shows.
\begin{example}\label{ex: chain sbm}(Chain SBM) We measure the clustering stability when applied to an SBM as described above.
	\begin{itemize}
		\item We use $r=8$ communities of size $100$ in the reduced model and use $\mu_1 I_2,\ldots,\mu_{r-1} I_2$ on the off-diagonal, where 
		\[ \mu_k=\frac{r-k}{r-1} \mu_1\]
		and $\mu_1 \in \{2,4,\ldots,100\}$ \bde(matrix \eqref{eq: reduced block model matrix}	shows such a model with $r=3$ and $\mu_1=20$)\ede. For small values of $\mu_1$, we expect a clustering into $r$ communities to be most stable, whereas for increasing $\mu_1$, the optimal number of clusters should decrease.  We compute the optimal number of clusters $k_\opt(\delta)$ provided by the SDA and the optimal number of clusters $k_\opt(g)$ as provided by the spectral gaps. Figure \ref{fig:kOpt} shows the results. We observe the expected behavior in both robustness measures, but the SDA tends to opt for a lower number of clusters for smaller values of $\mu_1$. Figure \ref{fig:kVals} shows the measures $\delta_6,\delta_7,\delta_8$ and $g_6,g_7,g_8$ for different values of $\mu_1$. As we expect, $\delta_8$ and $g_8$ are decreasing, $\delta_7$ and $g_7$ are increasing up to a certain point and then decreasing again.
		\item To compare the behavior of the above reduced model to the full SBM, we compute the same values for a SBM with $r=6$ communities of size $30$, edge probabilities $p_k = \mu_k/100$, where $\mu_1 \in [0,50]$. Figure \ref{fig:kOpt_SBM} shows the resulting optimal number of clusters, Figure \ref{fig:kVals_SBM} shows the distances to ambiguity and the spectral gaps for different values of $p_1$. We see how the results are affected by randomness in the graph generation but conclude that the behavior of both models is similar. 
	\end{itemize}
\end{example}
\begin{figure}
	\centering
	\include{SBM1}
	\caption{Optimal number of clusters in Example~\ref{ex: chain sbm} provided by the SDA and by the spectral gaps.}\label{fig:kOpt_SBM}
\end{figure}
It is remarkable that the optimal number of communities provided by both stability indicators does not differ by more than one. This suggests that the spectral gap is a reasonable stability indicator in this case, even though it does not take the special structure of weight matrices into account. We observe a similar behavior in different stochastic block models with a similar structure.

\subsection{Randomly distributed numbers around centers}
In the previous section we have compared the spectral gap to the distance to ambiguity as a stability measure. We have seen a similar behavior yet we cannot conclude which measure works better or worse. We now present an example where we can measure a success rate of the stability indicators.

Consider a given set of centers $m_1,\ldots,m_r$ with the intention to create clusters $C_1,\ldots,C_r$ around them. We generate a set of $n$ random numbers $x_i$ such that for each $i \in \{1,\ldots,n\}$, we first choose an index $j \in \{1,\ldots,r\}$ randomly (uniformly distributed) and then the random number $x_i$ which is normally distributed with expectation value $m_j$ and variance $1$. If the centers $m_1,\ldots,m_r$ are well separated, we expect that the numbers $x_1,\ldots,x_n$ are naturally clustered into $r$ communities.

For a graph with weights $w_{ij}$ related to the distances $|x_i-x_j|$, we then compute the SDA $\delta_k$ and the spectral gap $g_k$ for $k \in \{ r_{\min},r_{\max} \}$. Since the data set is constructed such that it should naturally have $r$ communities, we expect that $\arg\max_k \delta_k = r$ in most cases.

\begin{figure}
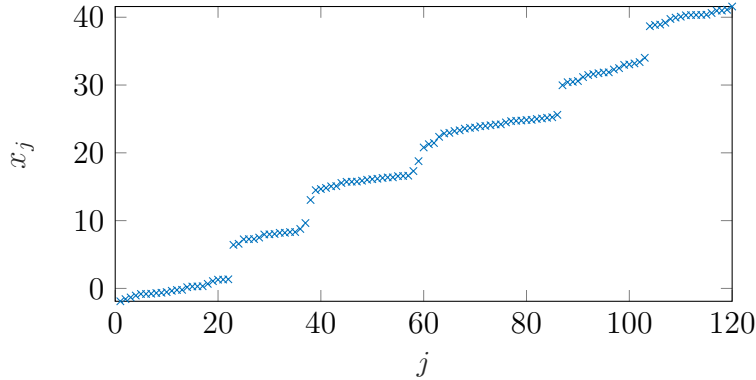

	\centering
	\include{randomsample}
	\caption{Random sample of $120$ numbers in six communities as  in Example~\ref{ex: random numbers}.}
	\label{fig:random number sample}
\end{figure}
\begin{example}\label{ex: random numbers}
	We generate $250$ samples of $n=120$ random numbers and $r=6$ groups. For each $i=1,\ldots,n$, we first choose a group $j \in \{1,\ldots,r\}$ randomly and then generate a random number $x_i$, normally distributed around $m_j$, where
	\[ m=\begin{pmatrix} 0 & 8 & 16 & 24 & 32 & 40 \end{pmatrix} \,. \]
	Figure \ref{fig:random number sample} shows one sorted random sample $(x_1,\ldots,x_n)$.
	
	In order to represent this data set by a graph, we set, following \cite{shi2000normalized}
	\[ f(x_i,x_j) = \exp(-\alpha(x_i-x_j)^2) \,, \]
	and then use
	\[ w_{ij} = \begin{cases}
	f(x_i,x_j)\,,\quad&\text{if }f(x_i,x_j)\ge \mathrm{tol} \,,\\
	0 \quad&\text{otherwise.}
	\end{cases} \]
	\bde The tolerance value (here $\mathrm{tol}=10^{-4}$) is used to neglect small weights and avoid a complete graph.\ede We denote by $g_{\opt}$ and $\delta_{\opt}$ the optimal number of clusters obtained by taking the spectral gaps and  the SDA, respectively, as stability indicators in the criterion of maximal stability. We used $r_{\min}=4$ and $r_{\max}=8$ as bounds for the number of communities. Table~\ref{tab:optimal values} reports the frequency of occurrence of $k$ optimal clusters for $4 \le k \le 8$ with both stability indicators for $\alpha = \frac{1}{2}$, and in Table~\ref{tab:optimal values 2} for $\alpha = \frac{1}{4}$. We conclude that the success of recognizing the number of communities strongly depends on the graph representation of the data set, but in both cases the SDA outperforms the spectral gap in finding the correct number of clusters.
\end{example}

\begin{table}[h!]
	\begin{center}
	\begin{tabular}{l | c | c | c|c|c}
		$k$ & $4$ & $5$ & $6$ & $7$ & $8$ \\
		\hline
		$g_{\opt}$ & 0.0\% & 0.0\% & \textbf{69.6\%} & 20.4\% & 0.0\% \\
		$\delta_{\opt}$  & 0.0\% & 0.0\% & \textbf{78.0\%} & 5.2\% & 6.8\% \\
		\hline
	\end{tabular}
	\caption{Frequency of $k$ being the optimal number of clusters provided by the spectral gaps and the SDA in Example~\ref{ex: random numbers} ($\alpha=1/2$).}\label{tab:optimal values}
	\end{center}
\end{table}
\begin{table}[h!]
	\centering
	\begin{tabular}{l | c | c | c|c|c}
		$k$ & $4$ & $5$ & $6$ & $7$ & $8$ \\
		\hline
		$g_{\opt}$ & 0.0\% & 0.0\% & \textbf{89.6\%} & 7.6\% & 2.8\% \\
		$\delta_{\opt}$  & 0.0\% & 0.0\% & \textbf{95.6\%} & 3.2\% & 1.2\% \\
		\hline
	\end{tabular}
	\caption{Frequency of $k$ being the optimal number of clusters provided by the spectral gaps and the SDA in Example~\ref{ex: random numbers} ($\alpha=1/4$).}\label{tab:optimal values 2}
\end{table}

\section{Conclusion} \label{sec:conclusion}
Stability of spectral clustering for undirected weighted graphs is measured by the {\it structured distance to ambiguity} (SDA) that is introduced and  studied here. We compare the SDA with the  spectral gap as a commonly accepted heuristic stability indicator. The spectral gap is characterized as an {\it unstructured} distance to ambiguity (Theorem~\ref{thm:sg}). In contrast to the spectral gap, the SDA takes into account the underlying graph structure, or more precisely, the $k$th SDA equals the distance of the given graph Laplacian to the set of Laplacian matrices of the same graph with perturbed weights and with coalescing $k$th and $(k+1)$-st eigenvalues. We derive a two-level algorithm to compute the SDA.

Stability indicators such as the structured or unstructured distance to ambiguity can be used for selecting the appropriate number of clusters, by maximizing the stability indicator (choosing the ``most stable'' clustering).
Our numerical experiments show that the spectral gap, although it ignores the underlying graph structure, often yields clustering results similar to the SDA. Tested on the stochastic block model, the optimal number of clusters delivered by both stability indicators usually differs by at most 1. We also present a particular example where the SDA outperforms the spectral gap in finding the optimal number of clusters.

The at most slightly better results of the SDA for choosing the number of clusters come along with a much higher computational cost. The purpose of this paper is thus not to replace the spectral gap by the SDA for general clustering computations.
The objective here was to provide a sound theoretical concept and an additional computational tool, against which clustering results obtained with computationally cheaper heuristic criteria can be tested in special situations. 
%
%
%
%
%

\section*{Acknowledgments}
We thank two anonymous Referees for their valuable comments and suggestions.

We thank Ulrike von Luxburg and Matthias Hein for a helpful and encouraging discussion.

Part of this work was developed during some visits to Gran Sasso Science Institute in L'Aquila and to the University of T\"ubingen.
The authors thank both institutions for the very kind hospitality.

N. Guglielmi thanks the INdAM GNCS (Gruppo Nazionale di Calcolo Scientifico) and the Italian M.I.U.R. (PRIN 2017) for financial support.

\bibliography{psgraphbib,biblioTesi}
\bibliographystyle{plain}

\end{document}

%% file: reducedSBM1.tex
%
%
\definecolor{mycolor1}{rgb}{0.00000,0.44700,0.74100}%
\definecolor{mycolor2}{rgb}{0.85000,0.32500,0.09800}%
\begin{tikzpicture}

\begin{axis}[%
width=3.229in,
height=1.434in,
at={(0.542in,0.442in)},
scale only axis,
xmin=2,
xmax=100,
xlabel style={font=\color{white!15!black}},
xlabel={$\mu{}_\text{1}$},
ymin=4,
ymax=9,
ytick={4, 5, 6, 7, 8, 9},
ylabel style={font=\color{white!15!black}},
ylabel={Optimal number},
axis background/.style={fill=white},
title style={font=\bfseries},
legend style={legend cell align=left, align=left, draw=white!15!black}
]
\addplot [color=mycolor1, draw=none, mark=x, mark options={solid, mycolor1}]
  table[row sep=crcr]{%
2	8\\
4	8\\
6	8\\
8	8\\
10	8\\
12	8\\
14	8\\
16	8\\
18	8\\
20	8\\
22	8\\
24	8\\
26	8\\
28	8\\
30	8\\
32	8\\
34	8\\
36	8\\
38	8\\
40	8\\
42	8\\
44	8\\
46	8\\
48	8\\
50	8\\
52	7\\
54	7\\
56	7\\
58	7\\
60	7\\
62	7\\
64	7\\
66	7\\
68	7\\
70	7\\
72	7\\
74	6\\
76	6\\
78	6\\
80	6\\
82	6\\
84	6\\
86	6\\
88	6\\
90	6\\
92	6\\
94	6\\
96	6\\
98	6\\
100	6\\
};
\addlegendentry{$k_{\text{opt}}(g)$}

\addplot [color=mycolor2, draw=none, mark=o, mark options={solid, mycolor2}]
  table[row sep=crcr]{%
2	8\\
4	8\\
6	8\\
8	8\\
10	8\\
12	8\\
14	8\\
16	8\\
18	8\\
20	8\\
22	8\\
24	8\\
26	8\\
28	8\\
30	8\\
32	8\\
34	8\\
36	8\\
38	7\\
40	7\\
42	7\\
44	7\\
46	7\\
48	7\\
50	7\\
52	7\\
54	7\\
56	7\\
58	7\\
60	7\\
62	7\\
64	7\\
66	7\\
68	6\\
70	6\\
72	6\\
74	6\\
76	6\\
78	6\\
80	6\\
82	6\\
84	6\\
86	6\\
88	6\\
90	6\\
92	6\\
94	6\\
96	6\\
98	6\\
100	5\\
};
\addlegendentry{$k_{\text{opt}}(\delta)$}

\end{axis}
\end{tikzpicture}%

%% file: reducedSBM2.tex
%
%
\definecolor{mycolor1}{rgb}{0.00000,0.44700,0.74100}%
\definecolor{mycolor2}{rgb}{0.85000,0.32500,0.09800}%
\definecolor{mycolor3}{rgb}{0.92900,0.69400,0.12500}%
\begin{tikzpicture}

\begin{axis}[%
width=1.917in,
height=1.698in,
at={(0.745in,0.229in)},
scale only axis,
xmin=0,
xmax=100,
xlabel style={font=\color{white!15!black}},
xlabel={$\mu{}_\text{1}$},
ymin=0,
ymax=150,
axis background/.style={fill=white},
legend style={at={(0.03,0.5)}, anchor=west, legend cell align=left, align=left, draw=white!15!black}
]
\addplot [color=mycolor1, line width=2.0pt]
  table[row sep=crcr]{%
2	138.91417056405\\
4	136.275579698112\\
6	133.507047071945\\
8	130.607316017362\\
10	127.573511944651\\
12	124.400376832337\\
14	121.078320388039\\
16	115.207336263449\\
18	111.524711133636\\
20	107.746235124778\\
22	103.872808038212\\
24	99.9073308297959\\
26	95.8547787776483\\
28	91.721827299888\\
30	87.5160248288324\\
32	83.2447684446865\\
34	78.9144668936859\\
36	74.5301186107894\\
38	70.0952655630713\\
40	65.6121424675819\\
42	61.0818617162654\\
44	56.504552882312\\
46	51.8794199755978\\
48	47.2047022631763\\
50	42.4775154920767\\
52	37.6935269263437\\
54	32.8463614918801\\
56	27.9265138717169\\
58	22.9192381157376\\
60	17.8000676204191\\
62	12.5243363804465\\
64	7.00462795769741\\
66	1.13574446606806\\
68	26.2652963876545\\
70	27.037805658996\\
72	26.3187136453689\\
74	27.0496360462109\\
76	27.7807285250483\\
78	28.5117814382719\\
80	30.9003486119776\\
82	29.9740559475667\\
84	30.7052151906143\\
86	31.4360734072373\\
88	32.1671347120343\\
90	32.8985604302267\\
92	33.6293083712117\\
94	34.3603762341103\\
96	37.0804183458238\\
98	17.0960185285057\\
100	4.22811792358357\\
};
\addlegendentry{$\delta{}_\text{8}$}

\addplot [color=mycolor2, line width=2.0pt]
  table[row sep=crcr]{%
2	4.00602634075609\\
4	8.01180163728797\\
6	12.0177134699852\\
8	16.0236282756828\\
10	20.029450591679\\
12	24.0353270067028\\
14	28.0412514651215\\
16	32.0471520954439\\
18	36.0529862547536\\
20	40.0588760497413\\
22	44.0647979122228\\
24	48.0709125528798\\
26	52.0765040688431\\
28	56.0823706173292\\
30	60.0883283442373\\
32	64.0941960823854\\
34	68.0999933395933\\
36	72.1059298882094\\
38	76.1118714226786\\
40	80.1177464069886\\
42	84.1235321503967\\
44	88.1294107701539\\
46	92.1357184637547\\
48	96.1412902583358\\
50	100.147192015154\\
52	104.152895797064\\
54	108.158863541436\\
56	112.164772239584\\
58	116.170753255562\\
60	120.176396780466\\
62	124.182376161325\\
64	128.18816350088\\
66	132.194062714401\\
68	65.6289713988822\\
70	61.6390431705186\\
72	57.5621006835707\\
74	53.4092672287593\\
76	49.2210313527239\\
78	45.0207632069236\\
80	40.8091675428217\\
82	36.5837245557818\\
84	32.3407243157007\\
86	28.0750716702639\\
88	23.7796938501415\\
90	19.4446747742121\\
92	15.0565032257458\\
94	10.5995353537491\\
96	6.06499530987428\\
98	1.46402482313141\\
100	3.1796234405854\\
};
\addlegendentry{$\delta{}_{\text{7}}$}

\addplot [color=mycolor3, line width=2.0pt]
  table[row sep=crcr]{%
2	2.91654997216958\\
4	5.83305675364798\\
6	8.74959136092547\\
8	11.6661383104447\\
10	14.5826453369732\\
12	17.4991630811326\\
14	20.4156713593456\\
16	23.3322841775637\\
18	26.2487698618149\\
20	29.1652386797314\\
22	32.0819160475086\\
24	34.9983486154387\\
26	37.9148143761895\\
28	40.8313273400984\\
30	43.747960005592\\
32	46.6643852413394\\
34	49.5809101766733\\
36	52.4974118751587\\
38	55.4139632807506\\
40	58.3304716898338\\
42	61.2473956347793\\
44	64.1635108367278\\
46	67.0800287602509\\
48	69.9965894561319\\
50	72.9132156046224\\
52	75.8295928304184\\
54	78.7461610880137\\
56	81.6626547311722\\
58	84.57919727946\\
60	87.4957265722359\\
62	90.4122770509164\\
64	93.3287314830557\\
66	96.2452454786443\\
68	99.1618074769637\\
70	102.078319712535\\
72	104.994824081255\\
74	107.91138791131\\
76	110.827931618813\\
78	113.744483156622\\
80	116.66093695992\\
82	119.577400224967\\
84	122.493983959899\\
86	125.410467893238\\
88	128.32700630283\\
90	131.243532792956\\
92	134.160108007875\\
94	137.076555723237\\
96	139.993208312788\\
98	142.909651398806\\
100	72.7579223819863\\
};
\addlegendentry{$\delta{}_{\text{6}}$}

\end{axis}

\begin{axis}[%
width=1.917in,
height=1.698in,
at={(3.268in,0.229in)},
scale only axis,
xmin=0,
xmax=100,
xlabel style={font=\color{white!15!black}},
xlabel={$\mu{}_\text{1}$},
ymin=0,
ymax=150,
axis background/.style={fill=white},
legend style={legend cell align=left, align=left, draw=white!15!black}
]
\addplot [color=mycolor1, line width=2.0pt]
  table[row sep=crcr]{%
2	137.160093002566\\
4	132.898829767823\\
6	128.63756653308\\
8	124.376303298337\\
10	120.115040063594\\
12	115.853776828851\\
14	111.592513594108\\
16	107.331250359365\\
18	103.069987124622\\
20	98.8087238898785\\
22	94.5474606551353\\
24	90.2861974203922\\
26	86.0249341856491\\
28	81.763670950906\\
30	77.502407716163\\
32	73.2411444814198\\
34	68.9798812466766\\
36	64.7186180119335\\
38	60.4573547771906\\
40	56.1960915424474\\
42	51.9348283077043\\
44	47.6735650729611\\
46	43.412301838218\\
48	39.1510386034749\\
50	34.8897753687318\\
52	30.6285121339887\\
54	26.3672488992456\\
56	22.1059856645025\\
58	17.8447224297594\\
60	13.5834591950165\\
62	9.32219596027314\\
64	5.06093272553017\\
66	0.799669490786954\\
68	3.16679256340624\\
70	3.25993352115345\\
72	3.35307447890074\\
74	3.44621543664805\\
76	3.53935639439526\\
78	3.63249735214247\\
80	3.72563830988976\\
82	3.81877926763696\\
84	3.91192022538421\\
86	4.00506118313144\\
88	4.09820214087873\\
90	4.19134309862594\\
92	4.28448405637313\\
94	4.37762501412042\\
96	4.47076597186759\\
98	4.56390692961483\\
100	2.69567966352116\\
};
\addlegendentry{$g_8$}

\addplot [color=mycolor2, line width=2.0pt]
  table[row sep=crcr]{%
2	1.39360874552006\\
4	2.78721749104015\\
6	4.18082623656024\\
8	5.57443498208036\\
10	6.96804372760052\\
12	8.36165247312047\\
14	9.7552612186407\\
16	11.1488699641607\\
18	12.5424787096808\\
20	13.9360874552009\\
22	15.329696200721\\
24	16.7233049462411\\
26	18.1169136917612\\
28	19.5105224372813\\
30	20.9041311828014\\
32	22.2977399283215\\
34	23.6913486738417\\
36	25.0849574193618\\
38	26.4785661648817\\
40	27.8721749104019\\
42	29.2657836559219\\
44	30.6593924014422\\
46	32.0530011469623\\
48	33.4466098924823\\
50	34.8402186380024\\
52	36.2338273835224\\
54	37.6274361290425\\
56	39.0210448745627\\
58	40.4146536200828\\
60	41.8082623656028\\
62	43.201871111123\\
64	44.595479856643\\
66	45.9890886021632\\
68	43.9211036037271\\
70	41.0534491145041\\
72	38.185794625281\\
74	35.3181401360581\\
76	32.4504856468349\\
78	29.5828311576121\\
80	26.7151766683889\\
82	23.8475221791659\\
84	20.9798676899429\\
86	18.11221320072\\
88	15.244558711497\\
90	12.3769042222739\\
92	9.50924973305093\\
94	6.64159524382789\\
96	3.77394075460485\\
98	0.906286265381996\\
100	1.961368223841\\
};
\addlegendentry{$g_7$}

\addplot [color=mycolor3, line width=2.0pt]
  table[row sep=crcr]{%
2	0.964385798578889\\
4	1.92877159715774\\
6	2.89315739573663\\
8	3.85754319431547\\
10	4.82192899289438\\
12	5.78631479147323\\
14	6.75070059005206\\
16	7.71508638863102\\
18	8.67947218720988\\
20	9.64385798578873\\
22	10.6082437843676\\
24	11.5726295829465\\
26	12.5370153815253\\
28	13.5014011801042\\
30	14.4657869786831\\
32	15.430172777262\\
34	16.3945585758408\\
36	17.3589443744196\\
38	18.3233301729985\\
40	19.2877159715774\\
42	20.2521017701563\\
44	21.2164875687352\\
46	22.180873367314\\
48	23.145259165893\\
50	24.1096449644717\\
52	25.0740307630508\\
54	26.0384165616296\\
56	27.0028023602083\\
58	27.9671881587874\\
60	28.9315739573661\\
62	29.895959755945\\
64	30.8603455545239\\
66	31.8247313531027\\
68	32.7891171516815\\
70	33.7535029502605\\
72	34.7178887488394\\
74	35.6822745474183\\
76	36.6466603459972\\
78	37.611046144576\\
80	38.5754319431549\\
82	39.5398177417339\\
84	40.5042035403125\\
86	41.4685893388915\\
88	42.4329751374704\\
90	43.3973609360493\\
92	44.3617467346282\\
94	45.3261325332071\\
96	46.290518331786\\
98	47.2549041303647\\
100	46.2579217051026\\
};
\addlegendentry{$g_6$}

\end{axis}
\end{tikzpicture}%

%% file: SBM1.tex
%
%
\definecolor{mycolor1}{rgb}{0.00000,0.44700,0.74100}%
\definecolor{mycolor2}{rgb}{0.85000,0.32500,0.09800}%
\begin{tikzpicture}

\begin{axis}[%
width=3.229in,
height=1.485in,
at={(0.542in,0.442in)},
scale only axis,
xmin=0.02,
xmax=0.5,
xlabel style={font=\color{white!15!black}},
xlabel={$p_1$},
ymin=2,
ymax=7,
ytick={2, 3, 4, 5, 6, 7},
ylabel style={font=\color{white!15!black}},
ylabel={Optimal number},
axis background/.style={fill=white},
legend style={at={(0.03,0.03)}, anchor=south west, legend cell align=left, align=left, draw=white!15!black}
]
\addplot [color=mycolor1, draw=none, mark=x, mark options={solid, mycolor1}]
  table[row sep=crcr]{%
0.02	6\\
0.04	6\\
0.06	6\\
0.08	6\\
0.1	6\\
0.12	6\\
0.14	6\\
0.16	6\\
0.18	6\\
0.2	6\\
0.22	6\\
0.24	6\\
0.26	6\\
0.28	5\\
0.3	5\\
0.32	5\\
0.34	5\\
0.36	5\\
0.38	5\\
0.4	5\\
0.42	5\\
0.44	5\\
0.46	4\\
0.48	4\\
0.5	4\\
};
\addlegendentry{$k_{\text{opt}}(g)$}

\addplot [color=mycolor2, draw=none, mark=o, mark options={solid, mycolor2}]
  table[row sep=crcr]{%
0.02	6\\
0.04	6\\
0.06	6\\
0.08	6\\
0.1	6\\
0.12	6\\
0.14	6\\
0.16	5\\
0.18	5\\
0.2	5\\
0.22	5\\
0.24	5\\
0.26	5\\
0.28	5\\
0.3	5\\
0.32	5\\
0.34	5\\
0.36	4\\
0.38	4\\
0.4	4\\
0.42	4\\
0.44	4\\
0.46	4\\
0.48	4\\
0.5	3\\
};
\addlegendentry{$k_{\text{opt}}(\delta)$}

\end{axis}
\end{tikzpicture}%

%% file: SBM2.tex
%
%
\definecolor{mycolor1}{rgb}{0.00000,0.44700,0.74100}%
\definecolor{mycolor2}{rgb}{0.85000,0.32500,0.09800}%
\definecolor{mycolor3}{rgb}{0.92900,0.69400,0.12500}%
\begin{tikzpicture}

\begin{axis}[%
width=1.917in,
height=1.698in,
at={(0.745in,0.229in)},
scale only axis,
xmin=0,
xmax=0.5,
xlabel style={font=\color{white!15!black}},
xlabel={$p_1$},
ymin=0,
ymax=30,
axis background/.style={fill=white},
legend style={at={(0.03,0.03)}, anchor=south west, legend cell align=left, align=left, draw=white!15!black}
]
\addplot [color=mycolor1, line width=2.0pt]
  table[row sep=crcr]{%
0.02	27.6503647389915\\
0.04	26.0752332136176\\
0.06	24.3544475515357\\
0.08	22.9143798973449\\
0.1	19.7895459085639\\
0.12	19.7958573734209\\
0.14	17.4850667496031\\
0.16	15.087378832427\\
0.18	14.655568596875\\
0.2	12.8202267464776\\
0.22	10.9681423835877\\
0.24	9.99008832732585\\
0.26	8.28493525968478\\
0.28	6.90428337370017\\
0.3	5.27820556541437\\
0.32	3.95148350647203\\
0.34	3.14444322694265\\
0.36	1.41905239538307\\
0.38	0.873906592554773\\
0.4	0.135413298563656\\
0.42	3.76822190084106e-14\\
0.44	0\\
0.46	1.00485917355762e-14\\
0.48	3.01457752067285e-14\\
0.5	2.26093314050464e-14\\
};
\addlegendentry{$\delta{}_\text{6}$}

\addplot [color=mycolor2, line width=2.0pt]
  table[row sep=crcr]{%
0.02	6.1765032587199\\
0.04	8.72338401405246\\
0.06	8.61718127901409\\
0.08	12.4529788571221\\
0.1	14.4205512920188\\
0.12	13.3136392451413\\
0.14	14.9696487751604\\
0.16	17.047555684297\\
0.18	17.2944702021049\\
0.2	17.2007149984997\\
0.22	18.1629669050989\\
0.24	19.8656545595986\\
0.26	20.9068810947422\\
0.28	20.4156754241944\\
0.3	21.4254420968507\\
0.32	20.4343618072529\\
0.34	21.1478890180935\\
0.36	21.2589209359751\\
0.38	11.26813172854\\
0.4	9.86590834965314\\
0.42	9.6336361989154\\
0.44	8.6599059267026\\
0.46	7.40598167870059\\
0.48	5.97323274513517\\
0.5	5.67693970474513\\
};
\addlegendentry{$\delta{}_{\text{5}}$}

\addplot [color=mycolor3, line width=2.0pt]
  table[row sep=crcr]{%
0.02	5.42341229864278\\
0.04	6.86248221361192\\
0.06	9.29657669036795\\
0.08	9.6134554919897\\
0.1	10.3403954428245\\
0.12	12.2800817316751\\
0.14	12.6199968259031\\
0.16	12.9507312262977\\
0.18	13.7817370793977\\
0.2	16.6676989525142\\
0.22	17.6582682871068\\
0.24	17.0099133750681\\
0.26	17.59039205896\\
0.28	18.9194465253864\\
0.3	19.8208204237329\\
0.32	20.0516431496783\\
0.34	19.7982813166959\\
0.36	21.9719463256783\\
0.38	21.2247080050527\\
0.4	22.7352893843873\\
0.42	22.6948765624693\\
0.44	22.0737423982566\\
0.46	24.194594147945\\
0.48	24.8952324838297\\
0.5	22.0444100711978\\
};
\addlegendentry{$\delta{}_{\text{4}}$}

\end{axis}

\begin{axis}[%
width=1.917in,
height=1.698in,
at={(3.268in,0.229in)},
scale only axis,
xmin=0,
xmax=0.5,
xlabel style={font=\color{white!15!black}},
xlabel={$p_1$},
ymin=0,
ymax=25,
axis background/.style={fill=white},
legend style={legend cell align=left, align=left, draw=white!15!black}
]
\addplot [color=mycolor1, line width=2.0pt]
  table[row sep=crcr]{%
0.02	19.8571297969566\\
0.04	18.934604421728\\
0.06	17.8983043944754\\
0.08	16.4493721320233\\
0.1	14.9427478747863\\
0.12	14.9154695555819\\
0.14	13.4029865177285\\
0.16	11.7213721651889\\
0.18	11.3585235330857\\
0.2	9.99161162078904\\
0.22	8.65605848202661\\
0.24	7.91481602662814\\
0.26	6.21027449897775\\
0.28	5.49697787967232\\
0.3	4.27328918304967\\
0.32	3.21953659647186\\
0.34	2.5421240115014\\
0.36	1.15782725354022\\
0.38	0.708831126930042\\
0.4	0.110187922435682\\
0.42	3.76822190084106e-14\\
0.44	0\\
0.46	1.00485917355762e-14\\
0.48	3.01457752067285e-14\\
0.5	2.26093314050464e-14\\
};
\addlegendentry{$g_6$}

\addplot [color=mycolor2, line width=2.0pt]
  table[row sep=crcr]{%
0.02	0.559979585697015\\
0.04	1.02439023146345\\
0.06	1.04039194161448\\
0.08	2.21968887100625\\
0.1	3.01279120168676\\
0.12	2.28160464915022\\
0.14	3.18615487188961\\
0.16	4.20827835509025\\
0.18	4.15742815034603\\
0.2	4.02353971966892\\
0.22	4.39762762607401\\
0.24	5.16994139287166\\
0.26	6.1695932692853\\
0.28	6.34479686704098\\
0.3	6.78817207641962\\
0.32	7.06574324517716\\
0.34	7.79047425788906\\
0.36	7.73425369872101\\
0.38	7.97657039625308\\
0.4	7.49576524432948\\
0.42	7.45596594972703\\
0.44	6.70048383490095\\
0.46	5.79423045173446\\
0.48	4.68911417336736\\
0.5	4.45197824997854\\
};
\addlegendentry{$g_5$}

\addplot [color=mycolor3, line width=2.0pt]
  table[row sep=crcr]{%
0.02	0.342224785043164\\
0.04	0.492451869846177\\
0.06	1.1786187165697\\
0.08	1.05043210040026\\
0.1	1.45861967070256\\
0.12	1.7945130579496\\
0.14	2.01130066784089\\
0.16	2.14663137220987\\
0.18	2.37036850776005\\
0.2	3.12659190533009\\
0.22	3.46834723788487\\
0.24	3.67327773201459\\
0.26	3.78073508837336\\
0.28	4.53292793406094\\
0.3	4.76318066523344\\
0.32	4.63687311150472\\
0.34	4.64513602490501\\
0.36	5.71818519138029\\
0.38	5.44381381254663\\
0.4	5.88366207588605\\
0.42	6.09856657356406\\
0.44	6.2365897520242\\
0.46	7.26479247857321\\
0.48	7.56295852811365\\
0.5	7.20517494216986\\
};
\addlegendentry{$g_4$}

\end{axis}
\end{tikzpicture}%

%% file: randomsample.tex
%
%
\definecolor{mycolor1}{rgb}{0.00000,0.44700,0.74100}%
\begin{tikzpicture}

\begin{axis}[%
width=3.229in,
height=1.544in,
at={(0.542in,0.383in)},
scale only axis,
xmin=0,
xmax=120,
xlabel style={font=\color{white!15!black}},
xlabel={\footnotesize{$j$}},
ymin=-1.898612957571,
ymax=41.5805044634241,
ylabel style={font=\color{white!15!black}},
ylabel={$x_j$},
axis background/.style={fill=white},
legend style={legend cell align=left, align=left, draw=white!15!black}
]
\addplot [color=mycolor1, draw=none, mark=x, mark options={solid, mycolor1}]
  table[row sep=crcr]{%
1	-1.898612957571\\
2	-1.57690009509428\\
3	-1.33004373221792\\
4	-1.05106157433722\\
5	-0.820779217943535\\
6	-0.819082113751756\\
7	-0.777414143093447\\
8	-0.680812488950661\\
9	-0.639234781206748\\
10	-0.540157686294364\\
11	-0.349160383611943\\
12	-0.230388370299027\\
13	-0.207121593013642\\
14	0.180370053309375\\
15	0.282141710692261\\
16	0.301511414725869\\
17	0.32617977410058\\
18	0.67225222604845\\
19	1.07517560156695\\
20	1.29159882584943\\
21	1.30292765477038\\
22	1.34681761101243\\
23	6.42212112123933\\
24	6.58271353996123\\
25	7.23942584082891\\
26	7.27042497973434\\
27	7.27944393813827\\
28	7.49973117427718\\
29	7.95502469618015\\
30	7.98032415359214\\
31	8.03271229682559\\
32	8.19762843917659\\
33	8.20673432993758\\
34	8.31429964752816\\
35	8.3166404850798\\
36	8.77283798497763\\
37	9.64380672500582\\
38	13.0443979172012\\
39	14.4880687933457\\
40	14.6611105349036\\
41	14.8236985278193\\
42	15.0578997600505\\
43	15.1001755670972\\
44	15.5540176421577\\
45	15.7240910447264\\
46	15.73770194197\\
47	15.7592678605446\\
48	15.8912843521312\\
49	16.0295274837552\\
50	16.1143826791298\\
51	16.1788704653106\\
52	16.2916606768929\\
53	16.3727700228243\\
54	16.3816449423174\\
55	16.5606263701962\\
56	16.5855353226582\\
57	16.6306369043827\\
58	17.3087583155233\\
59	18.772310380571\\
60	20.7912945676271\\
61	21.2716123019125\\
62	21.4480899386503\\
63	22.3611307956673\\
64	22.8575825667745\\
65	22.9503471889971\\
66	23.2134037235246\\
67	23.3211328767986\\
68	23.5678021513105\\
69	23.6805668261979\\
70	23.739766893791\\
71	23.915861228052\\
72	23.9725845533694\\
73	24.0708558518289\\
74	24.1815276455523\\
75	24.1982840016733\\
76	24.4846808718652\\
77	24.6741016615735\\
78	24.7251800453564\\
79	24.7608666220596\\
80	24.8410109240142\\
81	24.8611335216291\\
82	24.9905629474383\\
83	25.0720690917042\\
84	25.1735999302395\\
85	25.2476962791036\\
86	25.6078731059011\\
87	29.9651464641797\\
88	30.4149167054737\\
89	30.4760324382156\\
90	30.5900779411722\\
91	31.1704830812513\\
92	31.4758866365012\\
93	31.6117012908425\\
94	31.744168832191\\
95	31.8544323924327\\
96	31.8778792191108\\
97	32.2915692303795\\
98	32.4951943004756\\
99	32.9737940589087\\
100	33.0394415592402\\
101	33.1896518489601\\
102	33.3875809828851\\
103	34.0165934012355\\
104	38.6733460212805\\
105	38.8675636173678\\
106	38.9212315910224\\
107	39.1885120638426\\
108	39.7437400814894\\
109	39.9418366902784\\
110	40.0819723690089\\
111	40.2848163707486\\
112	40.3239385822918\\
113	40.3355455159793\\
114	40.3514705312957\\
115	40.3739914306515\\
116	40.6233206965122\\
117	40.9663951953436\\
118	40.9832328396497\\
119	41.1173092349514\\
120	41.5805044634241\\
};
\end{axis}
\end{tikzpicture}%